\documentstyle[12pt]{article}
\newcommand{\bt}{\begin{tabular}}
\newcommand{\et}{\end{tabular}}
\newcommand{\ba}{\begin{array}}
\newcommand{\ea}{\end{array}}
\newcommand{\be}{\begin{equation}}
\newcommand{\ee}{\end{equation}}
\newcommand{\ben}{\begin{enumerate}}
\newcommand{\een}{\end{enumerate}}

\newtheorem{Thm}{Theorem}[section]
\newtheorem{Prop}[Thm]{Proposition}
\newtheorem{Lem}[Thm]{Lemma}
\newtheorem{Rem}[Thm]{Remark}
\newtheorem{Ex}[Thm]{Example}
\newcommand{\dowod}{\noindent{\bf Proof:} }
\newcommand{\qed}{\hspace{5.5in} Q.E.D. }

\newcommand{\End}{\mbox{\rm End}\, }
\newcommand{\Mor}{\mbox{\rm Mor}\, }
\newcommand{\Vol}{\mbox{\rm Vol}\, }
\newcommand{\Ad}{\mbox{\rm Ad}\, }
\newcommand{\id}{\mbox{\rm id}\, }
\newcommand{\Ldwa}{\,\stackrel{2}{\bigwedge}}
\newcommand{\Ltrzy}{\,\stackrel{3}{\bigwedge}}
\newcommand{\w}{{\!}\wedge{\!}}
\newcommand{\h}{\hspace*{.1mm}}
\newcommand{\bu}{\bullet}
\newcommand{\sca}{g}
\newcommand{\krop}{{\mbox{\large\bf $\cdot$}}}
\newcommand{\lel}{\left\langle}
\newcommand{\rr}{\right\rangle}

\font \msa=msam10 scaled \magstep2

\font \msb=msbm10 scaled \magstep1
\font \msm=msbm10 scaled \magstep0
\newcommand{\rtimes}{\mbox{\msb o}\,}
\newcommand{\contr}{\mbox{\msa y}\,}
\newcommand{\bR}{\mbox{\msb R} }
\newcommand{\bC}{\mbox{\msb C} }
\newcommand{\bCm}{\mbox{\msm C} }

\font \eul=eufm10 scaled \magstep2
\font \eu=eufm10 scaled \magstep0

\newcommand{\gotG}{\mbox{\eul g}}
\newcommand{\goG}{\mbox{\eu g}}
\newcommand{\gotH}{\mbox{\eul h}}
\newcommand{\goH}{\mbox{\eu h}}
\newcommand{\gotN}{\mbox{\eul n}}
\newcommand{\goN}{\mbox{\eu n}}

\newcommand{\ar}{\alpha }
\newcommand{\tar}{\widetilde{\alpha}}
\newcommand{\br}{\beta }
\newcommand{\gr}{\gamma }
\newcommand{\dr}{\delta }
\newcommand{\er}{\varepsilon }
\newcommand{\lr}{\lambda }
\newcommand{\Om}{\Omega }

\begin{document}

\title{\bf Poisson structures on the Poincar\'{e} group}
\author{{\bf S. Zakrzewski}  \\
\small{Department of Mathematical Methods in Physics,
University of Warsaw} \\ \small{Ho\.{z}a 74, 00-682 Warsaw, Poland} }
\date{}
\maketitle

\begin{abstract}
An introduction to inhomogeneous Poisson groups is given.
Poisson inhomogeneous $O(p,q)$ are shown to be coboundary, the
generalized classical Yang-Baxter equation having only
one-dimensional right hand side. Normal forms of the classical
$r$-matrices for the Poincar\'{e} group (inhomogeneous $O(1,3)$)
are calculated.
\end{abstract}

\section{Introduction}
In this paper we give the proofs of facts announced in our
previous article \cite{PPgr}, in which we have presented a list
of 23 normal forms of classical $r$-matrices on the Poincar\'{e}
group.

It is remarkable that the classification of Poisson Poincar\'{e}
groups turns out to be completely analogous to the classification
of quantum Poincar\'{e} groups given in \cite{qpoi}. Recall also
that the main motivation of these investigations is the
potential possibility of deforming the relativistic symmetry.

The paper is organized as follows. In Sect.~1 we recall basic
definitions and facts concerning Lie bialgebras, especially in
the context of semi-direct product Lie algebras. In Sect.~2 we
prove that inhomogeneous $o(p,q)$ Lie algebras have the following
interesting features of simple Lie algebras: all Lie bialgebra
structures are coboundary (all Poisson Lie structures are of the
$r$-matrix type) and the subspace of invariants in the third
antisymmetric tensor power of the Lie algebra (the right hand
side of the {\em generalized classical Yang-Baxter equation})
is only one-dimensional. We formulate a set of equations
determining the classical $r$-matrix and present some solutions.

In Sect.~3 we restrict ourselves to the case of the Poincar\'{e}
Lie algebra (inhomogeneous $o(1,3)$) and present the main
result: a table of solutions. Sections 4 and 5 are devoted to
the proofs.

All vector spaces and Lie algebras considered in this paper are
real and finite-dimensional.

\section{Preliminaries}
\subsection{Modules}

Let $\gotG$ be a Lie algebra and let $E$ be a vector space.
Recall that $E$ is a $\gotG$-module if a bilinear map
$$ \gotG \times E\ni (X, u)\mapsto Xu\in E$$
is given, such that $[X,Y]u=X(Yu)-Y(Xu)$ for $X,Y\in \gotG$,
$u\in E$. We denote by $E_{\goG}$ the subspace of invariant elements:
$$ E_{\goG}:= \{ u\in E : Xu=0\;\;\mbox{for}\; X\in \gotG\}.$$
A {\em morphism} from a $\gotG$-module $E_1$ to a
$\gotG$-module $E_2$ is a linear map
$T\colon E_1\to E_2$ such that $T(Xu)= XT(u)$ for $X\in \gotG$,
$u\in E_1$. The linear space of morphisms from $E_1$ to $E_2$ is
denoted by $\Mor _{\goG} (E_1,E_2)$.
We have also the well known alternative terminology:
\begin{eqnarray*}
\gotG\mbox{-modules} & = & \mbox{representations of}\; \; \gotG \\
\mbox{morphisms of}\;\; \gotG\mbox{-modules} & = &
\mbox{intertwiners}.
\end{eqnarray*}
The tensor product of modules (representations) is naturally
defined. An important example of the $\gotG$-module is $\gotG$
itself with the adjoint action: $XY:=[X,Y]$. For the purpose of
this paper, the most important $\gotG$-module will be $\Ldwa
\gotG$.

\subsection{Cocycles and coboundaries}
Let $E$ be a $\gotG$-module. Linear map $f$ from $\gotG$ to $E$
is a {\em cocycle} ({\em on} $\gotG$ {\em with values in} $E$) if
$$ f([X,Y])=Xf(Y)-Yf(X)$$
for $X,Y\in \gotG$. The space of cocycles on $\gotG$ with values
in $E$ is denoted by $Z(\gotG , E)$. If $r\in E$, then the
linear map
$$ \gotG \ni X\mapsto (\partial r)(X):= Xr\in E$$
is said to be the {\em coboundary} of $r$. Coboundaries of
elements of $E$ form a subspace in $Z(\gotG ,E)$ which is
denoted by $B(\gotG ,E)$. We set
$$ H(\gotG ,E):= {Z(\gotG ,E)}/ {B(\gotG ,E)}.$$
Note that $H(\gotG ,E)=\{0\}$ if and only if each cocycle is a
coboundary.
The well known {\em Whitehead's lemma} says that for semisimple
$\gotG$ and arbitrary $\gotG$-module $E$ we have $H(\gotG
,E)=\{0\}$.

In order to approach the case of semi-direct product Lie
algebras, let us note the following useful (very simple) facts.
\ben
\item The restriction of a cocycle to a Lie subalgebra is a
cocycle (on this subalgebra).
\item If the restriction of a cocycle $\dr \in Z(\gotG ,E)$ to a
Lie subalgebra $\gotH$ is a coboundary, i.e. there exists $r\in
E$ such that $\dr (X)=Xr$ for $X\in \gotH$, then
\be\label{dr0}
\dr _0 := \dr - \partial r\in Z(\gotG ,E)\qquad
\mbox{satisfies}\;\;\; \dr _0 |_{\goH } = 0.
\ee
\item Let $\gotG = \gotN \rtimes \gotH$ (semidirect product;
$\gotN$ is the ideal) and let $\dr _0\colon\gotG \to E$ be a
linear map. Then
\be\label{ZMor}
\dr _0 \in Z(\gotG ,E),\;\;\dr_0|_{\goH} =0
\Longleftrightarrow
\dr _0|_{\goN}\in Z(\gotN ,E)\cap \Mor _{\goH}(\gotN ,E).
\ee
\item For $\gotG = \gotN\rtimes \gotH$ and $E:=\Ldwa\gotG$, we have
\be\label{E}
E = \Ldwa \gotN \oplus (\gotN \otimes \gotH )\oplus \Ldwa \gotH
\qquad (\mbox{$\gotH$-invariant decomposition})
\ee
\be\label{MorE}
\Mor _{\goH}(\gotN ,E)=\Mor _{\goH}(\gotN ,\Ldwa\gotN )\oplus
\Mor _{\goH}(\gotN , \gotN\otimes \gotH )\oplus \Mor
_{\goH}(\gotN ,\Ldwa \gotH ).
\ee
\een
\begin{Ex} \label{Ex1}
If \ $\gotH$ is semi-simple and $\gotG := \bR \oplus
\gotH$, then $H(\gotG ,\Ldwa \gotG ) = \{ 0\}$.
\end{Ex}
\dowod
If  $\dr \in Z(\gotG ,\Ldwa \gotG )$ then $\dr |_{\goH}\in
B(\gotH ,\Ldwa \gotG)$ (Whitehead's lemma), i.e. there exists
$r\in \Ldwa\gotG $ such that $\dr (X) = Xr$ for $X\in \gotH $.
Setting $\dr _0 := \dr -\partial r$ and using points 2 and 3
above we see that $\dr _0\in \Mor _{\goH}(\bR ,\Ldwa \gotG )$.
Note that
\be\label{zero}
\gotH _{\goH} =\{ 0\},\qquad (\Ldwa\gotH )_{\goH} = \{ 0 \}
\ee
(for the last equality, see e.g.~\cite{ccc}, Thm.I, p.~189). It
follows that
\be\label{R+h}
\Mor _{\goH}(\bR ,\Ldwa\gotG ) =\Mor _{\goH}(\bR , \gotH
)\oplus \Mor _{\goH}(\bR ,\Ldwa\gotH ) = \gotH _{\goH}\oplus
(\Ldwa\gotH )_{\goH} = \{ 0 \} .
\ee
\qed

Analogous fact holds of course for complex Lie algebras (with
$\bR $ replaced by $\bC$).

\subsection{Lie bialgebras}

Recall \cite{D:ham,D} that a {\em Lie bialgebra} is a Lie
algebra $\gotG$ together with a cocycle $\dr \colon \gotG\to
\Ldwa\gotG$ such that the dual map $\dr ^*\colon \Ldwa\gotG ^*
\to \gotG ^*$  is a Lie bracket on $\gotG^*$ (the dual of
$\gotG$). There is a 1--1 correspondence between Lie bialgebras
and connected, simply connected Poisson Lie groups
\cite{D:ham,D,S-T-S,Lu-We}.

A Lie bialgebra $(\gotG ,\dr )$ is said to be {\em coboundary}
if $\dr$ is a coboundary: $\dr =\partial r$, $r\in \Ldwa\gotG$.
Of course, Lie bialgebras $(\gotG ,\dr )$ with $\gotG$
semisimple are always coboundary. A non-semisimple Lie algebra
with the same property is provided by Example~\ref{Ex1} and the
following special case of it.
\begin{Ex}
Any Lie bialgebra structure on $\gotG = gl (n)$ is coboundary.
\end{Ex}
If $\dr =\partial r$ then $\dr ^*$ is a Lie
bracket if and only if
\be\label{YB}
[r,r]\in (\Ltrzy\gotG) _{\goG}
\ee
(the bracket used here is the Schouten bracket). Condition
(\ref{YB}) is called {\em generalized classical Yang-Baxter
equation} and $r$ is said to be a {\em classical $r$-matrix}.
The Lie bracket defined by $\partial r$ on $\gotG ^*$ equals
\be\label{nawr}
[\ar ,\br ]_r = r(\ar )\br - r(\br )\ar,\qquad \ar,\br \in \gotG
^*,
\ee
where $r(\cdot )\colon \gotG ^*\to \gotG $ is the contraction
with $r$ from the left:
\be
 r(\ar ) := \ar \contr r,\qquad \left\langle r(\ar ),\br
\right\rangle = \left\langle r,\ar
\otimes\br\right\rangle,\qquad \ar,\br \in \gotG ^*.
\ee
If $r$ is {\em triangular}, i.e. $[r,r]=0$, then
$r(\cdot )$ is a Lie algebra homomorphism:
\be\label{homr}
r([\ar ,\br ]_r)=[r(\ar ),r(\br )] .
\ee
This is a consequence of the following useful formula:
\be\label{formula}
\frac12 \left\langle [r,r],\ar\wedge \br\wedge \gr\right\rangle =
\left\langle [r(\ar ),r (\br )] - r([\ar ,\br ]_r), \gr\right\rangle
,\qquad \ar,\br\gr\in \gotG ^* .
\ee

\section{Inhomogeneous $o(p,q)$ algebras}

We consider a $(p+q)$-dimensional real vector space $V\cong \bR
^{p+q}$, equipped with a scalar product $\sca$ of signature
$(p,q)$. Let $\gotH := o(p,q)$ denote the Lie algebra of the
group $O(p,q)$ of endomorphisms of $V$ preserving $\sca$. Let
$\gotG := V\rtimes \gotH$ be the corresponding `inhomogeneous'
Lie algebra.

We recall that $\Ldwa V$ is naturally isomorphic to $\gotH$ as
a $\gotH$-module. The isomorphism is given by $\Om := \id
\otimes \sca $ (here $\sca $ is interpreted as a map from $V$ to
$V^*$). For $x,y\in V$ we set
\be\label{Omega}
 \Om _{x,y}:= \Om (x\wedge y) = x\otimes \sca (y) -
y\otimes \sca (x)\in \gotH\subset\End V .
\ee
 When working with a basis $e_1,\ldots ,e_{p+q}$ of $V$, we
shall use also the following notation
\be\label{Om}
\Om _{j,k} := \Om _{e_j ,e_k}=
 (\sca _{kl}e_j -\sca _{jl}e_k)\otimes e^l , \qquad j,k=1,\ldots
,p+q
\ee
(summation convention), where $e^1,\ldots , e^{p+q}$ is the dual
basis and $\sca _{jk}:= \sca (e_j,e_k)$.

\begin{Thm} \label{1}
For $\dim V >2$ we have
\be\label{H=0}
 H(\gotG ,\Ldwa \gotG )=\{ 0\}\qquad\qquad\mbox{({\em i.e. any
cocycle $\dr\colon \gotG\to \Ldwa \gotG$ is a coboundary})}
\ee
\be\label{inj}
 (\Ldwa\gotG )_{\goG} =\{ 0\}\qquad\qquad
\mbox{({\em  $r\mapsto  \partial r $ is injective}).}\hspace*{3.5cm}
\ee
\end{Thm}
\dowod Note that $\gotH =o(p,q)$ is semisimple for $p+q>2$ (it
is even simple, except the case $o(4,0)=o(3,0)\oplus o(3,0)$ and
$o(2,2)= o(2,1)\oplus o(2,1)$).
We first prove (\ref{inj}). Indeed, using (\ref{E})
with $\gotN = V$ and
\be\label{inva}
 (\Ldwa V)_{\goH }\cong (\gotH )_{\goH} = \{ 0\},\qquad (\Ldwa\gotH
)_{\goH}=\{ 0\}
\ee
(as in (\ref{zero})), we have
\be\label{VH}
(\Ldwa\gotG )_{\goH} = (V\otimes \gotH) _{\goH }
\cong \Mor _{\goH} (V,\gotH ).
\ee
The latter space is $\{ 0\}$ if $\dim V >3$ (because $V$ and
$\gotH$ are irreducible $\gotH$-modules of different dimension;
only $\gotH =o(4,0)$, $\gotH =o(2,2)$ are reducible, but in this
case the irreducible $\gotH $-submodules are of
dimension $3$). If $\dim V=3$, we have $\gotH \cong V$ and
$$\Mor _{\goH} (V,\gotH )= \Mor _{\goH} (V,V) \cong \bR .$$
It is easy to check that in this case the non-zero
$\gotH$-invariant element
\be\label{s}
s:=\er ^{jkl}e_j\otimes \Om _{kl}
\ee
 of $V\otimes \gotH\subset \Ldwa\gotG $
is not $V$-invariant (here $\er ^{jkl}$ is the usual antisymmetric
symbol). Concluding, $$(\Ldwa\gotG )_{\goG} =\{ 0\}.$$

To prove (\ref{H=0}), it is sufficient (in view
of (\ref{dr0}),(\ref{ZMor}) and semisimplicity of $\gotH$) to show
that if $\dr _0\in Z(V,\Ldwa\gotG )\cap \Mor _{\goH}
(V,\Ldwa\gotG )$ then $\dr _0 =\partial r$ for some $r\in (
\Ldwa\gotG)_{\goH}$. We shall show it first for $p+q>3$ (we
shall actually show that $\dr _0=0$). In this
case,
\be\label{zero1}
\Mor _{\goH} (V,\Ldwa V )=\Mor _{\goH} (V,\gotH )=\{ 0\}
\ee
(cf. remark after (\ref{VH})). We have also
\be\label{zero2}
\Mor _{\goH} (V,\Ldwa \gotH )=\{ 0\},
\ee
as a consequence of the following lemma.
\begin{Lem} Let $o(N)$ denote the orthogonal complex Lie algebra
acting in $\bC ^N$. Then
$$\Mor _{o(N)}(\bC ^N , \Ldwa o(N))=\{ 0\}\qquad \mbox{\em for}\
\ \  N\neq 3.$$
\end{Lem}
\dowod
It is sufficient to consider $N>3$. We consider two cases.
 \ben
  \item \fbox{$N=2n$}.\ \ Recall that weights of a
$\gotG$-module are functions on a basis of a Cartan subalgebra
of $\gotG$. For the $o(N)$-module $\bC ^N$, these functions are
non-zero at exactly one point. The basis may be chosen in such a
way that the non-zero values of these functions are $\pm 1$.

The weights of $o(N)\cong \Ldwa \bC ^N$ are sums of two
different weights of $\bC ^N$, hence either they are zero or
they are non-zero at exactly two points, where they have value
$\pm 1$ (note that this shows that $\Mor _{o(2n)}(\bC
^{2n},o(2n))=\{ 0\}$.)

The weights of $\Ldwa o(N)$ are sums of two different weights of
$o(N)$. The only weights having one-point support have values
$\pm 2$, hence there is no nontrivial intertwiner from $\bC ^N$
to $\Ldwa o(N)$.
 \item \fbox{$N=2n+1$}.\ \ For a suitably chosen $e_0\in \bC
^N$, we may identify $o(2n)$ as the subalgebra of $o(2n+1)$
stabilizing $e_0$, and acting on its orthogonal complement,
identified with $\bC ^{2n}$ (we can also choose the quadratic
form conveniently, if needed). We have
\be
\Mor _{o(2n+1)}(\bC ^{2n+1},\Ldwa (\Ldwa \bC ^{2n+1} ))\subset
\Mor _{o(2n)} (\bC\oplus\bC ^{2n},\Ldwa (\Ldwa \bC ^{2n+1} )).
\ee
Since
$$ \Ldwa (\bC\oplus\bC ^{2n} )\cong \bC ^{2n}\oplus\Ldwa \bC ^{2n}
 ,$$
$$ \Ldwa (\Ldwa (\bC\oplus\bC ^{2n}))\cong \Ldwa \bC ^{2n}
 \oplus (\bC ^{2n}\otimes \Ldwa \bC
^{2n} )\oplus \Ldwa (\Ldwa \bC ^{2n} ),$$
and
$$ \Mor _{o(2n)} (\bC , \Ldwa \bC ^{2n}) =
 \Mor _{o(2n)} (\bC , o(2n)) = (o(2n))_{o(2n)} = \{ 0\},$$
$$ \Mor _{o(2n)}(\bC , \bC ^{2n}\otimes o(2n))\cong \Mor
_{o(2n)} (\bC ^{2n},o(2n))=\{ 0\},$$
$$ \Mor _{o(2n)} (\bC , \Ldwa o(2n))\cong (\Ldwa
o(2n))_{o(2n)}=\{ 0\},$$
we have
$$ \Mor _{o(2n)} (\bC\oplus\bC ^{2n} ,\Ldwa o(N))\cong \Mor
_{o(2n)} (\bC ^{2n},\Ldwa o(N)).$$
If $T\colon \bC ^N\to \Ldwa o(N)$ is a nonzero $o(N)$-morphism,
then  $f|_{\bCm} =0$ and $ f|_{\bCm ^{2n}}$ is injective (since
$\bC ^{2n}$ is $o(2n)$-irreducible). Let $X$ be any element
of $o(N)$ which applied to $e_0$ gives a non-zero element of
$\bC ^{2n}$ (for instance $X=\Om _{j0}$). We obtain the contradiction
$$ 0\neq T(X e_0) =X T(e_0) =X (0)=0,$$
showing that $T$ has to be zero.
\een
\qed

\noindent
Due to (\ref{zero1}), (\ref{zero2}) and (\ref{MorE}),
\be
\Mor _{\goH} (V,\Ldwa \gotG )=\Mor _{\goH} (V,V\otimes\gotH ).
\ee
Using the fact that $\gotH$-modules $V$ and $\gotH$ are
isomorphic to their duals, we have
\be\label{kota}
 \Mor _{\goH} (V,V\otimes\gotH )\cong \Mor _{\goH} (\gotH
,V\otimes V )\cong \Mor _{\goH } (\gotH ,\Ldwa V ).
\ee
Here in the last equality we have used the following simple fact
(which may be easily proved using e.g.~\cite{Bourb},
\S~7,~Prop.~10).
\begin{Lem}
${\displaystyle \Mor _{o(N)} (o(N), \bC ^N \otimes _{\rm symm}\bC
^N ) = \{ 0\}\qquad }$ for $N>2$.
\end{Lem}
(Here the subscript `symm' refers to the symmetric part). It
follows that
\be\label{klam}
\Mor _{\goH} (V,\Ldwa \gotG )\cong
\Mor _{\goH} (\gotH ,\gotH )\cong \left\{ \ba{ll} \bR ^2 &
\;\;\mbox{if}\;\; p+q=4,\\
\bR & \;\;\mbox{otherwise.}
\ea\right.
\ee
The identity of $\gotH$ defines the following element $F_0$ of
$\Mor _{\goH} (V,V\otimes\gotH )$:
\be\label{F0}
V\ni x\mapsto F_0 (x):= \sca ^{jk}e_j \otimes \Om _{x,e_k}\in
V\otimes\gotH
\ee
($\sca ^{jk}$ is the contravariant metric). When $p+q=4$, the
Hodge star operation $* \colon \gotH\to \gotH $ given by
\be
*\Om _{x,z} := \sca (x)\wedge \sca (z) \contr \Vol,\qquad
(*\Om _{x,z})y= =\sca (x)\wedge\sca (y)\wedge\sca (z)\contr\Vol
=x\times y\times z
\ee
($\Vol $ is the volume element, $\times $ denotes the {\em
vector product} of three vectors) intertwines $\gotH $ with itself
and is not proportional to the identity. It defines another,
linearly independent from $F_0$, intertwiner from $V$ to
$V\otimes \gotH$:
\be
 F_1:= (\id \otimes *)\circ F_0.
\ee
Note that an element $F\in \Mor _{\goH} (V,\Ldwa \gotG )$
belongs to $Z(V,\Ldwa \gotG )$ if and only if the map
\be\label{symm}
V\times V\ni (x,y)\mapsto xF(y)\in \Ldwa \gotG
\ee
is symmetric. It is easily checked that $xF_0(y)$ is antisymmetric:
\be\label{xF0}
 xF_0(y) =  y\wedge x .
 \ee
If $p+q\neq 4$, it means that
\be\label{cap}
Z(V,\Ldwa \gotG )\cap \Mor _{\goH}(V,\Ldwa \gotG )=\{ 0\}.
\ee
 If $p+q=4$, one can show that
$$ xF_1(y)= \sca ^{jk} e_j \wedge (x\times y\times e_k),$$
which is also antisymmetric and linearly independent from
(\ref{xF0}). This shows that (\ref{cap}) holds also in this case.

Now let us consider the case $p+q=3$. Since $V\cong \Ldwa V\cong
\gotH$, we have
\be
\Mor _{\goH} (V,\Ldwa V ) \cong \bR ,\qquad  \Mor _{\goH}
(V,V\otimes \gotH ) \cong \bR ,\qquad \Mor _{\goH} (V,\Ldwa
\gotH ) \cong \bR
\ee
(cf. also (\ref{kota})). Note that the symmetry
of (\ref{symm}) for $F\in \Mor _{\goH} (V,\Ldwa\gotG )$
means the symmetry condition separately for each of its three
components  (in the decomposition (\ref{MorE}) with $\gotN =V$).
The first component is proportional to
\be\label{first}
V\ni x\mapsto T(x):=\sca (x)\contr \Vol\in \Ldwa V.
\ee
The symmetry is trivially satisfied in this case. The second
component, proportional to (\ref{F0}) satisfies the symmetry of
(\ref{symm}) if and only if it is zero, by (\ref{xF0}). The
third component is proportional to
\be\label{trud}
(\Om\otimes\Om )(T\otimes T)T.
\ee
One can show by a direct calculation, that
(\ref{trud}) does not satisfy the symmetry condition. We
conclude that in the case when $p+q=3$,
\be\label{cap1}
Z(V,\Ldwa \gotG )\cap \Mor _{\goH}(V,\Ldwa \gotG )=\{ \bR \cdot T\}.
\ee
But $T=-\frac12 \partial s$, where $s$ is given by (\ref{s}).
This ends the proof of the theorem.

\qed

In view of this theorem, the classification of Lie bialgebra
structures on $\gotG=V\rtimes \gotH$ consists in a description
of equivalence classes (modulo ${\rm Aut}\, \gotG$) of $r\in
\Ldwa\gotG$ such that $[r,r]\in (\Ltrzy \gotG )_{\goG }$.

Each $r\in \Ldwa \gotG$ has a decomposition
$$ r = a + b + c ,$$
corresponding to the  decomposition (\ref{E})
$$ \Ldwa \gotG = \Ldwa V \oplus (V\wedge \gotH ) \oplus \Ldwa
\gotH .$$
We have also the following decomposition of the Schouten bracket
\be\label{Sch}
 [r,r] = 2[a,b] + (2[a,c] +[b,b]) +2[b,c] +[c,c],
\ee
corresponding to the decomposition
$$ \Ltrzy \gotG = \Ltrzy V \oplus (\Ldwa V\wedge \gotH ) \oplus
(V\wedge \Ldwa\gotH ) \oplus \Ltrzy \gotH .$$
Note that
\be\label{dec}
 (\Ltrzy \gotG )_{\goG}= (\Ltrzy V )_{\goG}\oplus
(\Ldwa V\wedge \gotH )_{\goG} \oplus
(V\wedge \Ldwa\gotH )_{\goG} \oplus (\Ltrzy \gotH )_
{\goG}.
\ee
We shall show that this space is one-dimensional for $p+q>3$.
Note that the isomorphism $\Om $ defines a canonical
$\gotH$-invariant element of $(\Ldwa V)^*\otimes \gotH $, or,
using the identification of $V$ and $V^*$, a canonical
$\gotH$-invariant element of $\Ldwa V \otimes \gotH$. We shall
denote this element again by $\Om $. It is given by
\be\label{Ome}
 \Om = \sca ^{jl} \sca ^{km} e_j\wedge e_k \otimes \Om _{l,m}
\ee
(in any basis). This element is also $V$-invariant:
$$ x\Om =-\sca ^{jl} \sca ^{km} e_j\wedge e_k \wedge (e_l
x_m-e_mx_l)=0 \qquad \mbox{for} \;\; x\in V .$$
  \begin{Thm}
 If $\dim V > 3$ then  $(\Ltrzy \gotG )_{\goG}=
(\Ldwa V\wedge \gotH )_{\goG} = \bR \cdot \Om $.
\end{Thm}
\dowod
We calculate all terms in (\ref{dec}).
\ben
\item If $w\in \Ltrzy V$ is $\gotH$-invariant, then
$$V\ni x\mapsto \sca (x)\contr w\in \Ldwa V$$
belongs to $\Mor _{\goH} (V,\Ldwa V )$. {}From (\ref{zero1}) it
follows that  $w=0$. Hence $(\Ltrzy V)_{\goG}=\{ 0\} $.
\item The second component in (\ref{dec}) is contained in
(\ref{klam}). We already know that $\Om$ is $\gotG$-invariant.
If $p+q=4$, the second (linearly independent) $\gotH$-invariant
element $(\id\otimes *)\Om $ of $\Ldwa V\otimes \gotH$ is not
$V$-invariant:
$$ x(\id \otimes * )\Om =
-\sca ^{jl} \sca ^{km} e_j\wedge e_k \wedge (*\Om _{lm})x
= e_j\wedge e_k \wedge (e^j\wedge e^k \wedge \sca (x)\contr \Vol
)= 2\sca (x)\contr \Vol .$$
It follows that $(\Ldwa V \wedge \gotH )_{\goG} =\bR \cdot \Om $.
\item The third component in (\ref{dec}) is zero by
(\ref{zero2}).
\item We shall show that $(\Ltrzy \gotH )_{V} =\{ 0\}$. If $w\in
\Ltrzy \gotH$ is $V$-invariant, then
$$ 0= x w \in V\wedge \Ldwa \gotH \qquad \mbox{for}\;\; x\in V ,$$
hence
$$ 0 = \xi \contr xw \qquad \mbox{for}\;\; x\in V,\; \xi\in V^* .$$
Since $\xi \contr xw = - \omega _{\xi ,x} \contr w$, where
$\omega _{\xi ,x}\in \gotH ^*$ is defined by $\omega _{\xi ,x}
(A):= \left\langle \xi ,A x\right\rangle $, we have
$$ \ar\contr w =0\qquad \mbox{for}\;\; \ar\in \gotH ^*$$
(elements of the form $\omega _{\xi ,x}$ span $\gotH ^*$),
hence $w=0$.
\een
\qed

{}From this result and (\ref{Sch}) it follows that Lie bialgebra
structures on $\gotG$ are (for $p+q> 3$) in one-to-one
correspondence with $r=a+b+c\in \Ldwa \gotG$ such that
\begin{eqnarray}
  {} [c,c] & = & 0 \label{cc} \\
  {} [b,c] & = & 0 \label{bc} \\
  {} 2[a,c] + [b,b] & = & t\h \Om \qquad (t\in \bR ) \label{bb} \\
  {} [a,b] & = & 0 . \label{ab}
\end{eqnarray}
Equation (\ref{cc}) means that $c$ is a {\em triangular}
$r$-matrix on $\gotH$ (this is the semi-classical counterpart of
a known theorem \cite{inho} excluding the case when the
homogeneous part $H$ is $q$-deformed). Equation (\ref{bc}) tells
that $b$, as a map from $\gotH ^*$ to $V$, is a cocycle:
\be\label{bcoc}
b([\ar ,\br ]_c)= c(\ar )b(\br )-c(\br )b(\ar )\qquad \mbox{for}
\;\ar ,\br \in \gotH ^*,
\ee
the Lie bracket on $\gotH ^*$ being defined by the triangular
$c\in \Ldwa \gotH$ as in (\ref{nawr}):
\be\label{bra}
 [\ar ,\br ]_c = c(\ar )\br - c(\br )\ar
\ee
and the action of $\gotH ^*$ on $V$ is defined using the
homomorphism from $\gotH ^*$ to $\gotH$ given by $c$:
$$ c(\ar ):= \ar \contr c \in \gotH\qquad \mbox{for}\;\; \ar\in
\gotH ^*$$
(as in (\ref{homr})). To get (\ref{bcoc}) one can use
(\ref{formula}) with $\ar,\br \in \gotH ^*$, $\gr\in V^*$.

Here are  some particular solutions of (\ref{cc})--(\ref{ab}).
\ben
\item $a=0$, $b=0$, $c\in\Ldwa \gotH$ triangular.
 \item $b=0$, $c=0$, $a\in \Ldwa V$ arbitrary. This type of
solutions we call `soft deformations' \cite{abel}.
 \item $a=0$, $c=0$, $[b,b]=t\Om $. There is a family of
solutions of the latter equation, parameterized by vectors in
$V$. Namely, for each $x\in V$,
\be\label{bx}
 b_x:= F_0(x)= \sca ^{jk} e_j\otimes \Om _{x,e_k} = \frac12
g(x)\contr \Om
\ee
satisfies this equation with $t=-\sca (x,x)$ ($F_0$ is defined
in (\ref{F0})). Moreover, since $[x,b_x]=0$ (easy calculation),
\be\label{bX}
 b=b_x + x\wedge X ,\qquad X\in \gotH _x \;\;(\mbox{stabilizer
of $x$ in $\gotH$})
\ee
satisfies $[b,b]=[b_x,b_x] = -\sca (x,x)\Om$. Indeed, $[x\w
X,x\w X]=0$ and
$$[x\w X ,b_x]=x\w Xb_x - X\w [x,b_x] = x\w b_{Xx}=0.$$
\een
Note the following two properties of $b$ given in (\ref{bX}) for
$x\neq0$:
\be\label{X-2}
[a,b]=0\;\;\; \Longleftrightarrow \;\;\; (X-2)a\in x\wedge V\qquad
\mbox{for}\;\; a\in \Ldwa V
\ee
\be\label{X-1}
(-v)b = x\wedge (X-1)v\qquad\qquad \mbox{for}\;\; v\in V
\ee
(the first follows from $[a,b] = x\wedge (X-2)a$).
\begin{Prop} \label{noa}
Suppose $b$ is given by (\ref{bX}) with $x\neq 0$. If $X$ has no
eigenvalue 1 on $V$ and no eigenvalue 2 on $\Ldwa V$, then for
any solution $a$ of (\ref{ab}), $r=a+b$ can be transformed to $b$
by a suitable internal automorphism of $\gotG$.
\end{Prop}
\dowod
Since $X$ preserves $x\wedge V$ and $(X-2)$ is invertible on
$\Ldwa V$, the right hand side of (\ref{X-2}) is equivalent to
$a\in x\wedge V$. Since $X-1$ is invertible, $(-v)b$ runs over
$x\wedge V$ when $v$ runs over $V$. It follows that $(\Ad
_{-v}\otimes \Ad _{-v})(a+b)= (a+(-v)b) + b$ is equal $b$
for some $v\in V$ (see also (\ref{cech})).

\qed

Of course, a generic $X$ will satisfy the assumptions of the
above proposition.

\section{The case of the Poincar\'{e} group}

We now fix $(p,q)=(1,3)$. It means that $V\cong \bR ^{1+3}$ is
the four-dimensional Minkowski space-time, $\gotH =o(1,3)\cong
sl (2,\bC ) $ is the Lorentz Lie algebra and $\gotG $ is the
Poincar\'{e} Lie algebra.

We are interested in classifying the solutions of
(\ref{cc})--(\ref{ab}) up to the automorphisms of $\gotG$. In
particular, $c$ can be always chosen to be a normal form of a
triangular classical $r$-matrix on the Lorentz Lie algebra, as
listed in \cite{repoi}.  In the next section, for each such
non-zero $c$, we shall solve (\ref{bc})--(\ref{ab}) completely
(up to an automorphism).  Moreover, we shall find all solutions
with $c=0$ provided $t=0$. The results are shown in
Table~\ref{tabela} below.  Let us explain the notation. We
introduce the standard generators of $\gotH = sl(2,\bC )$:
$$H=\frac12 \left[ \ba {cc} 1 & 0 \\ 0 & -1 \ea \right] ,\qquad
X_+=\left[ \ba {cc} 0 & 1 \\ 0 & 0 \ea \right] ,\qquad
X_-=\left[ \ba {cc} 0 & 0 \\ 1 & 0 \ea \right]. $$
The action of $X\in sl(2,\bC )$ on a vector $v\in
V$ is given by $X(v):= Xv + vX^+$, the space $V$ being
identified with the set of hermitian $2\times 2$ matrices, where
$X^+$ is the hermitian conjugate of $X$. We fix the Lorentz basis
$e_0,e_1,e_2,e_3$ in $V$ given by the standard Pauli matrices:
$$e_0=\left[\ba{cc} 1 & 0 \\ 0 & 1 \ea\right] ,\;\;\;\;\;
e_1=\sigma _1 =\left[\ba{cc} 0 & 1 \\ 1 & 0 \ea\right] ,\;\;\;\;\;
e_2 =\sigma _2 = \left[\ba{cc} 0 & -i \\ i & 0 \ea\right] ,\;\;\;\;\;
e_3=\sigma _3= \left[\ba{cc} 1 & 0 \\ 0 & -1 \ea\right] .$$
We denote by $J$ the multiplication by the imaginary unit in
$\gotH$. As acting on $V$, the basic generators of $\gotH$ are given by
\be\label{HJH}
H=L_3=\Om _{30} = e_0\otimes e^3 + e_3\otimes e^0,\qquad JH =
-M_3=\Om _{21}=e_1\otimes e^2 - e_2\otimes e^1
\ee
\be\label{X+}
X_+=\Om _{10}+\Om _{13}=\Om _{e_1,e_+},\qquad
JX_+=\Om _{02}+\Om _{32}=\Om _{e_+,e_2},
\ee
\be\label{X-}
X_-=\Om _{10}+\Om _{31}=\Om _{e_1,e_-},\qquad
JX_-=\Om _{20}+\Om _{32}=\Om _{e_2,e_-} .
\ee
It is also convenient to introduce the light-cone vectors
$e_{\pm } := e_0 \pm e_3$.

The table lists 21 cases labelled by the number {\sf N} in the
last column.  In the forth column (labelled by {\#}) we indicate
the number of essential parameters (more precisely -- the
maximal number of such parameters) involved in the deformation.
This number is in many cases less than the number of parameters
actually occurring in the table.  The final reduction of the
number of parameters can be achieved using two following
one-parameter groups of automorphisms of $\gotG$:
\ben
 \item the group of dilations: $(v,X)\mapsto (\lr v,X)$
 (in cases 1,2,3,4,6),
 \item the group of internal automorphisms generated by $H$ and
the group of dilations (in cases 11,12,15,17,18).
\een
(the table looks more concise before the final reduction).
\begin{Rem} The table below differs a little from the table
announced in~\cite{PPgr}. Some errors are corrected and the
presentation is improved. Solutions (\ref{be0+}), (\ref{be1+})
are now presented separately (they are not included in the
table: they form the known part of the not yet solved problem
$[b,b]=t\Om ,\;t\neq 0$) and are now supplemented by
(\ref{be1+a}).
\end{Rem}

\pagebreak

\begin{table}[t]
\begin{center}
\bt{|@{\h}c@{\h}|@{\h}c@{\h}|@{\h}c@{\h}|@{\h}c@{\h}|@{\h}c@{\h}|}
  \hline
 $c$ & $b$ & $a$  & \# & {\sf N} \\
  \hline\hline
 $\gr JH\w H$ & $0$ & $\ar e_+\w e_- + \tar e_1\w e_2 $ & 2 & 1 \\
  \hline
 $JX_+\w X_+ $ & $\br_1 b_{e_+}+\br _2 e_+\w   JH $ & 0 &  1 & 2 \\
   \cline{2-5}
 & $\br b_{e_+}$ & $\ar e_+ \w e_1 $ & 1 & 3 \\
   \cline{2-5}
 & $\br (e_1\w X_+ + e_2\w JX_+)$ & $e_+ \w
 (\ar _1 e_1 + \ar _2 e_2)- \br ^2 e_1\w e_2 $ &  2  & 4 \\
 \hline
  $H\w X_+ - $  &  &  &    &  \\
 $ JH\w JX_+ +  $  & $0$ & $0$ &  1  & 5 \\
  $ \gr JX_+\w X_+ $  & & &   &  \\
 \hline
 $H\w X_+$ & $\br_1 b_{e_2}+\br _2 e_2\w X_+$ & 0 & 1  & 6 \\
 \hline
$0$ & $b_{e_+} + \br e_+\w JH$ & $0$ & 1 & 7 \\
   \cline{2-5}
  & $ b_{e_+} + \br e_+\w X_+$ & $ 0$ & 1  & 8 \\
    \cline{2-5}
 & $ e_1\w (X_+ +\br JX_+) + $\hspace{8mm}
 & $ \ar e_+\w e_2$ &  2  & 9 \\
 & $+ e_+\w (H+\sigma X_+),\;\;\; $ $\sigma =0,\pm 1$ & &  &  \\
    \cline{2-5}
 & $ e_1\w JX_+ + e_+\w X_+$ & $ \ar _1 e_-\w e_1 + \ar _2 e_+\w
e_2$ & 2 & 10 \\
   \cline{2-5}
& $ e_2\w X_+$ & $ \ar _1e_+\w e_1 + \ar _2e_-\w e_2$ &  1
 & 11 \\
 \cline{2-5}
 & $ e_+\w X_+$ & $e_-\w (\ar e_+ +\ar_1e_1 +\ar_2e_2)
 +\tar e_+\w e_2$ & 3 & 12 \\
  \cline{2-5}
& $ e_0\w JH$ & $ \ar_1e_0\w e_3 + \ar _2 e_1\w e_2$ &  2  & 13 \\
   \cline{2-5}
 \cline{2-5}
& $ e_3\w JH$ & $ \ar_1e_0\w e_3 + \ar _2 e_1\w e_2 $ &  2  & 14 \\
 \cline{2-5}
 & $ e_+\w JH$  &
  $\ar_1e_0\w e_3 + \ar _2 e_1\w e_2  $ &  1 & 15 \\
      \cline{2-5}
& $ e_1\w H$ & $ \ar_1e_0\w e_3 + \ar _2e_1\w e_2 $ & 2  & 16 \\
      \cline{2-5}
& $ e_+\w H$ & $\ar e_1\w e_2 + \ar _1 e_+\w e_1 $ & 1 & 17 \\
      \cline{2-5}
 & $ e_+\w (H+\br JH)$ &
  $\ar e_1\w e_2 $ & 1
  & 18 \\
    \cline{2-5}
 & $0$ & $e_1\w e_+$ &  0  & 19 \\
      \cline{3-5}
 &  & $e_1\w e_2$ &  0 & 20 \\
      \cline{3-5}
 &  & $ e_0\w e_3 + \ar e_1\w e_2$ &  1 & 21 \\
 \hline
\et
\caption{Normal forms of $r$ for $c\neq 0$ or $t=0$.}
\label{tabela}
\end{center}
\end{table}

In the case when $c=0$ and $t\neq 0$, the only solutions
we know are based on formula (\ref{bX}). We describe them now.
We shall use yet another standard generators of $\gotH$:
$$M_i = \er _{ijk} e_k\otimes e{^j},\qquad  L_i = e_0\otimes
e{^i} + e_i\otimes e{^0}\qquad (i,j,k=1,2,3).$$
 If we set $x:= e_0$ in (\ref{bx}), we obtain
 $$ b_{e_0}  =  e_1\wedge L_1 + e_2\wedge L_2 + e_3\wedge L_3 ,$$
 which is the known \cite{luki} classical $r$-matrix
corresponding to so called $\kappa$-deformation. More generally,
using (\ref{bX}), we have
\be\label{be0+}
 b = b_{e_0} +\lr e_0\wedge M_3
\ee
(any element of $\gotH _x\cong o(0,3)$ can be rotated to $\lr
M_3$).  Since $M_3$ has only imaginary eigenvalues, adding $a$
we do not obtain essentially different solutions,
cf.~Prop.~\ref{noa}.

Taking $x=e_1$ in (\ref{bx}), we obtain another solution
\be\label{be1}
 b_{e_1} = e_0\wedge L_1 -  e_2\wedge M_3 + e_3\wedge M_2
\ee
(this one is $M_1,L_2,L_3$-invariant). There are three types of
elements in $\gotH _x\cong o(1,2)$, according to the sign of the
Killing form. We have thus three types of perturbations
(\ref{bX}) of (\ref{be1}):
\be\label{be1+}
 b=b_{e_1} +\lr e_1\wedge Y,\qquad Y=M_1\;\;\mbox{or}
\;\;Y=M_1+L_3\;\; \mbox{or}\;\;Y= L_3=H .
\ee
In the first two cases, adding $a$ does not yield new solutions,
since $M_1$ has only imaginary eigenvalues and $M_1+L_3$ is
nilpotent. Since non-zero eigenvalues of $H$ are $\pm 1 $, adding
$a$ in the third case we can obtain a nontrivial modification
when $\lr =\pm 1,\pm 2$. We obtain then the following four
families of solutions:
\be\label{be1+a}
b= b_{e_1} \pm k e_1\wedge H + \ar  e_{k}\wedge
e_{\pm },\qquad k=1,2
\ee
(using the automorphisms generated by $H$, we can assume that
$\ar =\pm 1$).

\section{The proof for $c\neq 0$}

The four types of non-zero triangular $c$ in the table above
are\hspace{-.3mm}
taken from \cite{repoi}. We consider each case separately.
{\hspace{-.2mm}}We denote by $(H^*,JH^*,X_{\pm}^*,JX_{\pm }^*)$
the basis dual to $(H,JH,X_{\pm},JX_{\pm})$.

\subsection{$c=JH\wedge H$}

First we calculate brackets (\ref{bra}) of basis elements and
write down corresponding cocycle condition (\ref{bcoc}). We do
not consider pairs of elements from the subset
$\{X_+^*,JX_+^*,X_-^*,JX_-^*\}$, since for them the corresponding
condition (\ref{bcoc}) is trivial. We have
$$\ba{rclrcl}
{}[JH^*,H^*]_c & = & 0  & 0 & = & Hb(H^*)+JH b(JH^*) \\
{}[X_{\pm}^*,H^*]_c & = & \pm JX_{\pm}^* & b(JX_{\pm}^* ) & = &
\pm JHb(X_{\pm}^*) \\
{}[JX_{\pm}^*,H^*]_c & = & \mp X_{\pm}^* & b(X_{\pm}^*) & = &
 \mp JHb(JX_{\pm}^*) \\
{}[X_{\pm}^*,JH^*]_c & = & \pm X_{\pm}^* & b( X_{\pm}^*) & = &
\mp  Hb(X_{\pm}^*)\\
{}[JX_{\pm}^*,JH^*]_c & = & \pm JX_{\pm}^* &\hspace{3cm} b(
JX_{\pm}^*) & = &  \mp  Hb(JX_{\pm}^*) .
\ea
$$
 Due to (\ref{HJH}), $H (V)\cap JH (V)=\{ 0\}$, hence
the last four formulas imply
$$ b(X_{\pm}^*)=0,\qquad b(JX_{\pm}^*)=0.$$
For the same reason, the first equation,
$$ Hb(H^*)=-JH b(JH^*),$$
has the obvious solution $b(H^*)\in \ker H$, $b(JH^*)\in \ker
JH$, which can be written as follows
$$ b(H^*) = JH (v),\qquad b(JH^*)= -H(v),\qquad v\in V,$$
or
$$ b = H\wedge JH (v) - JH \wedge H(v) = vc.$$
Using the internal automorphism $\Ad _{-v}=\id - v$ of $\gotG$,
we can transform $r=a+b+c$ into
\be\label{cech}
(\Ad _{-v}\otimes \Ad _{-v}) r = r + (-v)r + (v\otimes v)r =
(a + (-v) b + (v\otimes v)c) + (b + (-v)c) +c,
\ee
hence we can always set $b=0$.

The last equation to solve is (\ref{bb}) with $b=0$. Since $\Om$
represents the isomorphism of $\Ldwa V$ and $\gotH$, it has rank
equal $6 =\dim\gotH$ (as an element of the tensor product of
$\Ldwa V$ and $\gotH$). The Schouten bracket $[c,a]$ is a tensor
of rank at most 2, because
$$ [JH\wedge H,a] = JH\wedge [H,a]-H\wedge [JH,a ].$$
It follows that $t=0$, hence equation (\ref{bb}) reduces to
$[c,a]=0$, i.e.
$$ [H,a]=0,\qquad [JH,a]=0.$$
It is clear that $a\in \Ldwa V$, considered as an element of
$\gotH$ has to be a combination of $H$ and $JH$. Going back to
$\Ldwa V$ (and using (\ref{HJH})), we obtain
$$ a = \lr e_0\wedge e_3 +\mu e_1\wedge e_2,\qquad \lr ,\mu \in
\bR .$$
Since we can multiply our solution $r=c+a$ by any number, we
obtain the first case of the table.

It is easy to check that $X\in\gotH$ and $Xc\equiv [X,c]=0$
implies that $X$ is a combination of $H$ and $JH$.
Such $X$ gives rise to a group of internal automorphisms of
$\gotG$, leaving $c$ invariant. These automorphisms leave
invariant also $a$, hence they cannot be used to a further
reduction of $a$.

\subsection{$c= JX_+\wedge X_+$}

First we calculate the brackets (formula (\ref{bra})) of basis
elements (only those contributing to the cocycle condition):
\begin{eqnarray}
{}[JX_+^*,X_+^* ]_c & = & 2H^* \\
{} [JX_+^*,H^* ]_c & = & -2X_-^* \\
{} [X_+^*,H^* ]_c & = & -2JX_-^* \\
{} [X_+^*,JH^* ]_c & = & 2X_-^* \\
{} [JX_+^*,JH^* ]_c & = & -2JX_-^*
\end{eqnarray}
and $X_-^*$, $JX_-^*$ are central elements. It follows that the
cocycle condition (\ref{bcoc}) reads
\begin{eqnarray}
2b(H^*) & = & X_+b(X_+^*) + JX_+b(JX_+^*) \label{HXX}\\
2b(X_-^*) & = & -X_+b(H^*)  \\
2b(JX_-^*) & = & JX_+b(H^*)  \\
2b(X_-^*) & = & -JX_+b(JH^*)  \label{XJH} \\
2b(JX_-^*) & = & -X_+b(JH^*)  \label{JXJH}
\end{eqnarray}
and $b(X_-^*)$, $b(JX_-^*)$ are $X_+$- and $JX_+$-invariant (the
latter property is already a consequence of
(\ref{XJH})-(\ref{JXJH}), since $X_+\circ JX_+=0=JX_+\circ
X_+$). We recall (cf. (\ref{X+})) that
\be
  X_+ x=2x^-e_1 + x^1e_+,\;\;\;\; JX_+x= -2x^-e_2
-x^2e_+\;\;\;\;\mbox{for}\;\; x=x^+e_++x^-e_- +x^1e_1 +x^2e_2.
\ee
To solve (\ref{HXX})--(\ref{JXJH}) we can just set $b(X_+^*)=x$,
$b(JX_+^*)=y$, where $x,y\in V$ are arbitrary vectors and then
\begin{eqnarray}
b(H^*) & = & \frac12 ( X_+x +JX_+y)= x^-e_1 -y^-e_2 +\frac12
(x^1 -y^2)e_+ \\
b(X_-^*) & = & -\frac12 X_+b(H^*) = -\frac14 (X_+)^2x=-\frac12
x^-e_+ \\
b(JX_-^*) & = & \frac12 JX_+b(H^*) = \frac14 (JX_+)^2y=\frac12
y^-e_+ .
\end{eqnarray}
Equations (\ref{XJH})-(\ref{JXJH}) will be satisfied by
$b(JH^*)=:z$ if
\begin{eqnarray*}
x^-e_+ & = & JX_+z= -2z^-e_2 - z^2e_+ \\
-y^-e_+ & = & X_+z = 2z^-e_1+ z^1e_+ ,
\end{eqnarray*}
i.e. $z = z^+e_+ - y^-e_1 - x^-e_2$ with arbitrary $z^+\in \bR $.
We have thus solved (\ref{bc}) completely (the solution is
parameterized by $x,y\in V$ and $z^+\in \bR $).

Now we are going to solve (\ref{bb}). Using formula
(\ref{Ome}) with the basis $e_+, e_-, e_1, e_2$, we have
\be\label{Omeg}
\Om = e_-\wedge e_+ \otimes H - 2e_1\wedge e_2\otimes JH +
e_-\wedge e_1 \otimes X_+ + e_2\wedge e_-\otimes JX_+ +
e_+\wedge e_1\otimes X_- + e_+\wedge e_2\otimes JX_- .
\ee
We shall compute terms on the left hand side of (\ref{bb}) which are
proportional to $e_-\wedge e_1\otimes X_+$, $e_-\wedge e_2\otimes
X_+$, $e_2\wedge e_-\otimes JX_+$. Note that they may come only
from $[b,b]$. Indeed,  $[X_+,a]$ and $[JX_+,a]$ are combinations
of $e_+\wedge e_1$, $e_+\wedge e_2$, $e_+\wedge e_-$, $e_1\wedge
e_2$, while
\be\label{ca}
[c,a]= JX_+\wedge [X_+,a] - X_+\wedge [JX_+,a] .
\ee
Using the general form of $b$,
$$ b=H\wedge b(H^*) + JH\wedge b(JH^*)+
X_+\wedge b(X_+^*) + JX_+\wedge b(JX_+^*) +
X_-\wedge b(X_-^*) + JX_-\wedge b(JX_-^*),$$
it is clear that the terms in $[b,b]$ which contain $X_+$ are
the following:
$$ 2\left( [H,X_+]\wedge b(H^*)\wedge b(X_+^*) + [JX_+,JH]\wedge
b(JX_+^*)\wedge b(JH^*)\right) +  2X_+\wedge [b(X_+^*),b] =$$
$$ = 2X_+\wedge \left( b(H^*)\wedge b(X_+^*) + b(JX_+^*)\wedge
b(JH^*) \right) + 2X_+\wedge \left( b(H^*)\wedge Hx  +  b(JH^*)\wedge
JHx + \right.
$$
$$ \left. +  b(X_+^*)\wedge X_+x + b(JX_+^*)\wedge
JX_+x + b(X_-^*)\wedge X_-x + b(JX_-^*)\wedge JX_-x \right) .$$
Now we substitute previously computed solutions, neglecting
terms which do not contribute to the factor at $e_-\wedge e_1$,
$e_-\wedge e_2$. We have (apart from $2X_+$)
$$ (x^-e_1 -y^-e_2)\wedge x + y\wedge (-y^-e_1 -x^-e_2) +
(x^-e_1 -y^-e_2)\wedge ( -x^-e_-) + x\wedge 2x^-e_1 +
y\wedge (-2x^-e_2) .$$
It is easy to write the part of $[b,b]$, proportional to
$2X_+\wedge e_-\wedge e_1$:
\be
2X_+\wedge e_-\wedge e_1\cdot (-{x^-}^2 - {y^-}^2 + {x^-}^2 +
2{x^-}^2)=
2X_+\wedge e_-\wedge e_1\cdot (2{x^-}^2 - {y^-}^2)
\ee
and to $2X_+\wedge e_-\wedge e_2$:
\be
2X_+\wedge e_-\wedge e_2\cdot (x^-y^- - y^-x^- - x^-y^- -2x^-y^-)=
2X_+\wedge e_-\wedge e_2\cdot (-3x^-y^-).
\ee
Similar calculation shows that the term proportional to
$2JX_+\wedge e_2\wedge e_-$ is
\be
2JX_+\wedge e_2\wedge e_-\cdot (2 {y^-}^2 - {x^-}^2).
\ee
Looking at (\ref{Omeg}), we see that
$$ 2{x^-}^2 - {y^-}^2 = 2 {y^-}^2 - {x^-}^2,\qquad x^-y^-=0,$$
which means that $x^- =0=y^-$ and $t=0$. In particular,
$b(X_-^*)=0=b(JX_-^*)$ and
\be\label{b1}
 b=X_+\wedge (x^+e_+ + x^1e_1+x^2e_2) +
 JX_+\wedge (y^+e_+ +y^1e_1+y^2e_2) +
 H\wedge \frac12 (x^1 -y^2)e_+  + JH\wedge z^+e_+ .
\ee
We shall simplify this general form, using appropriate
automorphisms of $\gotG $. First, note that
\be\label{-vc}
 (-v)c= JX_+ v\wedge X_+ + JX_+\wedge X_+v= X_+\wedge
(2v^-e_2 + v^2e_+)  + JX_+\wedge (2v^-e_1 + v^1e_+),
\ee
hence transforming $b$ into $b+(-v)c$ (as in (\ref{cech})) we
may assume that $x^+=0=y^+$ and $x^2+y^1 =0$ in (\ref{b1}).
Secondly, the one-parameter group of internal automorphisms
generated by $JH$ leaves $c$ invariant and transforms $b$
according to $\dot{b}=JH\cdot b$, i.e.
$$ \dot{x}^1=x^2 -y^1=\dot{y}^2,\qquad
\dot{x}^2=-(x^1+y^2)=-\dot{y}^1,$$
or
$$(x^1-y^2)^{\krop}=2(x^2+y^1),\qquad (x^2-y^1)^{\krop}=
 -2(x^1+y^2)$$
and $(x^1+y^2)^{\krop}=0=(x^2-y^1)^{\krop}$, hence one can afford
$x^2-y^1=0$. This implies $x^2=0=y^1$ (we already had
$x^2+y^1=0$) and we have the following simplified form of $b$:
\be\label{b2}
 b=x^1X_+\wedge e_1 +
 y^2JX_+\wedge e_2 +
 \frac12 (x^1 -y^2) H\wedge e_+  + z^+JH\wedge e_+ .
\ee
Now we can finally solve (\ref{bb}). We have
$$[b,b]=(x^1+y^2)\left( JX_+\wedge e_+\wedge (y^2e_2 +2z^+e_1)
-X_+\wedge e_+\wedge (x^1e_1 +2z^+e_2)\right) .$$
For $a=e_+\wedge (\ar _1e_1+\ar_2e_2)+e_-\wedge (\br_1e_1+\br_2e_2)+
\gr e_-\wedge e_+ + \dr e_1\wedge e_2$ we have also (see (\ref{ca}))
$$ [c,a]=JX_+\wedge (\br _1 e_-\wedge e_+ +2\br _2e_1\wedge e_2
-2\gr e_+\wedge e_1 +\dr e_+\wedge e_2) +
$$
$$-X_+\wedge (2\br _1 e_1\wedge e_2 -\br _2 e_-\wedge e_+ + 2\gr
e_+\wedge e_2 + \dr e_+\wedge e_1 ).
$$
It follows that $2[c,a]+[b,b]=0$ if and only if
$\br_1 = \br_2 = \gr = 0$ and
$$ z^+(x^1+y^2)=0,\qquad -2\dr = y^2(x^1+y^2)=x^1(x^1+y^2).$$
There are two possibilities:
\ben
\item $x^1+y^2=0$, $\dr =0$, i.e.
\be\label{2-3}
 b = x^1(X_+\wedge e_1 -JX_+\wedge e_2 + H\wedge e_+)  +
z^+JH\wedge e_+ , \qquad a=e_+\wedge (\ar _1e_1+\ar_2e_2),
\ee
\item $x^1+y^2\neq 0$, $z^+=0$, $x^1=y^2$, $\dr = -(x^1)^2$, i.e.
\be\label{case4}
 b = x^1(X_+\wedge e_1 +JX_+\wedge e_2), \qquad a=e_+\wedge (\ar
_1e_1+\ar_2e_2) - (x^1)^2e_1\wedge e_2.
\ee
\een
Of course, (\ref{case4}) is the case 4 in the table. Note, that
since $b$ in (\ref{2-3}) is $JH$-invariant (and $c$ also is),
one can transform $a$ to the following form:
$$ a = \ar e_+\wedge e_1 ,$$
because $JH$ generates rotations in the $e_1,e_2$-plane.
If $z^+=0$, we obtain case 3 in the table. If $z^+\neq 0$
we can get rid of $a$ in (\ref{2-3}) as follows. First we transform
the whole $r$ as in (\ref{cech}) with $v=v^1e_1+v^2e_2$, which
gives new $b$ (cf.~(\ref{-vc})) and $a$:
\begin{eqnarray}
b & = & x^1(X_+\wedge e_1 -JX_+\wedge e_2 + H\wedge e_+)  +
z^+JH\wedge e_+  + X_+\wedge v^2e_+  + JX_+\wedge v^1e_+ \\
a & = & e_+\wedge (\ar _1e_1+\ar_2e_2) + x^1(v^1e_+\wedge e_1
+ v^2e_+\wedge e_2)+z^+(v^2e_1-v^1e_2)\wedge e_+ .
\end{eqnarray}
We choose $v^1,v^2$ such that $a=0$. Now observe that
$X_+c=0$, $JX_+c=0$, $X_+b=-z^+JX_+\wedge e_+$ and
$JX_+b=z^+X_+\wedge e_+$, hence the automorphism groups
generated by $X_+$ and $JX_+$ change only $v^1$ and $v^2$,
respectively,  according to
$$ \frac{d}{dt} v^1 = -2z^+, \qquad \; \frac{d}{ds} v^2 = 2z^+$$
(parameters $t$ and $s$ correspond, respectively, to $X_+$ and
$JX_+$). Using these transformations we can afford $v^1=0=v^2$
(due to $z^+\neq 0$), which is the case 2 in the table.

\subsection{$c=H\wedge X_+ -JH\wedge JX_+ +\gr  JX_+\wedge X_+$}

Calculation of brackets (\ref{bra}) and corresponding cocycle
condition (\ref{bcoc}) gives
$$\ba{rclrcl}
{}[H^*,X_-^*] & = & 0 &  0 & = & X_+ b(X_-^*) \\
{}[H^*,JX_-^*] & = & 0 &  0 & = & X_+ b(JX_-^*) \\
{}[JH^*,X_-^*] & = & 0 &  0 & = & JX_+ b(X_-^*) \\
{}[JH^*,JX_-^*] & = & 0 &  0 & = & JX_+ b(JX_-^*) \\
{}[H^*,JH^*] & = & 0 & 0 & = & X_+ b(JH^*)+JX_+ b(H^*) \\
{}[X_+^*,X_-^*] & = & -X_-^* & \hspace{2cm}
 b(X_-^*) & = & (H+\gr JX_+)b(X_-^*) \\
{}[X_+^*,JX_-^*] & = & -JX_-^* & b(JX_-^*) & =
   & (H+\gr JX_+)b(JX_-^*) \\
{}[JX_+^*,X_-^*] & = & -JX_-^* & -b(JX_-^*) & =
   & (JH+\gr JX_+)b(X_-^*) \\
{}[JX_+^*,JX_-^*] & = & X_-^* & b(X_-^*) & = & (JH+\gr JX_+)b(JX_-^*)
\ea
$$
$$\ba{rclrcl}
{}\!\! [H^*,X_+^*] & \!\! =\!\!  & H^*+2\gr JX_-^* & \hspace{2mm}
 b(H^*) + 2\gr b(JX_-^*) & \!\! =\!\!
  & X_+b(X_+^*)+ (H+\gr JX_+)b(H^*) \\
{}\!\! [JH^*,X_+^*] & \!\! =\!\!  & JH^*-2\gr X_-^* &
 b(JH^*) - 2\gr b(X_-^*) & \!\! =\!\!  & \! -JX_+b(X_+^*)+
(H+\gr JX_+)b(JH^*) \\
{}\!\! [H^*,JX_+^*] & \!\! =\!\!  & JH^*+2\gr X_-^* &
 b(JH^*) + 2\gr b(X_-^*) & \!\! =\!\!
   & X_+b(JX_+^*)-(JH+\gr X_+)b(H^*) \\
{}\!\! [JX_+^*,JH^*] & \!\! =\!\!  & H^*-2\gr JX_-^* &
 b(H^*) - 2\gr b(JX_-^*) & \!\! =\!\!
   & JX_+b(JX_+^*)+(JH+\gr X_+)b(JH^*)
\ea
$$
$$
\ba{rclrcl}
{}[X_+^*,JX_+^*] & = & -2\gr H^* & \hspace{.7cm}
2\gr b(H^*) & = & (JH+\gr X_+)b(X_+^*)+ (H+\gr JX_+)b(JX_+^*)
\ea
$$
(and $[X_-^*,JX_-^*]=0$). We have suppressed the subscript `$c$'
in the bracket.
The first four equations imply that equations
from the sixth to the ninth take the form
\begin{eqnarray*}
 b(X_-^*) & = & Hb(X_-^*) \\
 b(JX_-^*) & = & Hb(JX_-^*) \\
 -b(JX_-^*) & = & JHb(X_-^*) \\
 b(X_-^*) & = & JHb(JX_-^*) .
\end{eqnarray*}
First two of the above equations imply that $b(X_-^*)$ and
$b(JX_-^*)$ are proportional to $e_+$ and then the last two
equations imply $b(X_-^*)=0=b(JX_-^*)$.
What remains is the following set of equations:
\begin{eqnarray}
X_+ b(JH^*) & = & -JX_+ b(H^*) \label{3.0}\\
 b(H^*) & = & X_+b(X_+^*)+ (H+\gr JX_+)b(H^*) \label{3.1}\\
 b(JH^*) & = & -JX_+b(X_+^*)+ (H+\gr JX_+)b(JH^*) \label{3.2}\\
 b(JH^*)  & = & X_+b(JX_+^*)-(JH+\gr X_+)b(H^*) \label{3.3}\\
 b(H^*)  & = & JX_+b(JX_+^*)+(JH+\gr X_+)b(JH^*) \label{3.4}\\
2\gr b(H^*) & = & (JH+\gr X_+)b(X_+^*)+ (H+\gr JX_+)b(JX_+^*)
\label{3.5} .
\end{eqnarray}
Knowing that $\ker X_+ =\left\langle e_+,e_2\right\rangle$,
$X_+e_1=e_+$, $X_+e_-=2e_1$,
$\ker JX_+ =\left\langle e_+,e_1\right\rangle$,
$JX_+e_2=-e_+$, $JX_+e_-=-2e_2$, one can easily solve (\ref{3.0}):
$$ b(H^*)=\ar e_+ + \br e_1 + \lr e_2,\qquad
b(JH^*)=\mu e_+ + \lr e_1 + \rho e_2, \qquad
\ar,\br,\lr,\mu,\rho\in \bR .$$
Setting $b(X_+^*)=: x$, we can write (\ref{3.1}) as follows:
$$ \ar e_+ + \br e_1 + \lr e_2 = x^1e_+ 2x^-e_1 + \ar e_+ - \lr
e_+ .$$
It means that $\lr = 0 = x^1$, $\br = 2x^-$. {}From (\ref{3.2}) we get
$$ \mu e_+ + \rho e_2 = x^2e_+ +2x^-e_2 + \mu e_+ - \gr \rho
e_+ ,$$
hence $\rho = 2x^-$, $x^2= \gr \rho = 2\gr x^-$. Recall that we
have now
$$ b(H^*)=\ar e_+ + 2x^- e_1 ,\qquad
b(JH^*)=\mu e_+ +  2x^- e_2.$$
Setting $b(JX_+^*)=:y$, we get from (\ref{3.3})
$$ \mu e_+ 2x^- e_2 = y^1 e_+ + 2y^- e_1 +2x^-e_2- 2\gr x^- e_+,$$
hence we get $y^-=0$ and $y^1 = \mu + 2\gr x^-$.  Equation
(\ref{3.4}) yields
$$ \ar e_+ + 2x^- e_1 = -y^2e_+ + 2x^-e_1 ,$$
hence $y^2 = -\ar $. Finally, (\ref{3.5}) yields
$$2\gr (\ar e_+ + 2x^- e_1 ) = x^2 e_1+2\gr x^-e_1 + y^+e_+ -
\gr y^2 e_+.$$
Since $x^2=2\gr x^-$, $y^2 =-\ar$, from this equation we get
$y^+=\gr\ar$. Concluding, the general solution of (\ref{bcoc}) is
$$
b= H\w (\ar e_+ + 2x^- e_1)+JH\w (\mu e_+ +  2x^- e_2) +
$$
$$+X_+\w (x^+e_+ + x^-e_- + 2\gr x^-e_2) + JX_+\w
(\gr\ar e_+ +(\mu + 2\gr x^-)e_1 -\ar e_2).
$$
Comparing this with $(-v)c$ for a general $v\in V$,
$$(-v)c = H\w (v^1 e_+ + 2v^- e_1)+JH\w (v^2 e_+ +  2v^- e_2) +
$$
$$ +X_+\w ((\gr v^2-v^+)e_+ + v^-e_- + 2\gr v^-e_2) + JX_+\w
(\gr v^1 e_+ +(v^2 + 2\gr v^-)e_1 -v^1 e_2),
$$
it is easy to see that $b=(-v)c$ for $v^1=\ar$, $v^2=\mu $,
$v^-=x^-$, $v^+=\gr\mu -x^+$. Therefore we can always assume
that $b=0$.

Now we shall show that $[c,a]=0\Longrightarrow a=0$ for $a\in
\Ldwa V$. Indeed,
$$ [c,a] = H\w X_+a - JH\w JX_+a -X_+\w (Ha+\gr JX_+a)+JX_+\w
(JHa +\gr X_+a)$$
is zero if and only if $X_+a=0$, $JX_+a=0$, $Ha=0$ and $JHa=0$.
But the commutant of $\{ H, X_+\}$ in $\gotH$ is zero.

\subsection{$c=H\wedge X_+ $}

We calculate brackets (\ref{bra}) relevant for the  cocycle condition
(\ref{bcoc}):
$$\ba{rcllrcl}
{}[H^*,X_+^*] & = & H^*         &1^{\bu}
              &  b(H^*) & =  & X_+b(X_+^*)+ Hb(H^*) \\
{}[H^*,JX_+^*] & = & JH^*
              &2^{\bu} &  b(JH^*) & =  & X_+b(JX_+^*) \\
{}[H^*,X_-^*] & = & 0           &3^{\bu} & 0 & = & X_+ b(X_-^*) \\
{}[H^*,JX_-^*] & = & 0          &4^{\bu} &  0 & = & X_+ b(JX_-^*) \\
{}[X_+^*,JH^*] & = & 0          &5^{\bu} & 0 & = & Hb(JH^*) \\
{}[H^*,JH^*] &=& -2JX_-^*       &6^{\bu}
              & 2b(JX_-^*) & = & -X_+ b(JH^*) \\
{}[X_+^*,JX_+^*] & = & JX_+^*  \hspace{2cm} &7^{\bu} &
\hspace{1cm}  b(JX_+^*) & = & -Hb(JX_+^*) \\
{}[X_+^*,JX_-^*] & = & -JX_-^*  &8^{\bu} & b(JX_-^*) & =
              & Hb(JX_-^*) \\
{}[X_+^*,X_-^*] & = & -X_-^*    &9^{\bu} & b(X_-^*) & = & Hb(X_-^*) .
\ea
$$
It follows from $7^{\bu}$-$9^{\bu}$ that $b(JX_+^*)= \ar e_-$,
$b(JX_-^*)= \br e_+$, $b(X_-^*)= \gr e_+$ and this implies
$3^{\bu}$-$4^{\bu}$. {}From $2^{\bu}$ we obtain $b(JH^*)=2\ar e_1$
which implies  $5^{\bu}$. $6^{\bu}$ means $2\br e_+=-X_+(2\ar
e_1)=-2\ar e_+$, hence $\br = -\ar$. The only remaining
condition is $1^{\bu}$:
$$ (1-H)b(H^*)=X_+ b(X_+^*).$$
Denoting $b(X_+^*)=:x$, $b(H^*)=:y$ we obtain
$$ 2y^-e_- +y^1e_1 + y^2 e_2 = x^1e_+ + 2x^- e_1,$$
i.e. $x^1=y^-=y^2 =0$, $ y^1=2x^-$. The general solution of the
cocycle condition is therefore
$$b=
H\w (y^+e_+ +2x^- e_1)+JH\w 2\ar e_1 + X_+\w (x^+e_+ + x^-e_- +
x^2e_2) + \ar JX_+\w e_- + \gr X_-\w e_+ - \ar JX_-\w e_+ .$$
Adding to this
$$ (-v)c=H\w X_+v - X_+\w Hv =
H\w (v^1e_+ +2v^- e_1) - X_+\w (v^+e_+ - v^-e_-)
$$
for a suitable $v\in V$, we get a simpler form of $b$:
\be
b=2\ar JH\w e_1 + x^2X_+\w e_2 + \ar JX_+\w e_- + \gr X_-\w e_+
- \ar JX_-\w e_+ .
\ee
We have $b= \ar b_0 + \br b_1 + \gr b_2$, where $\ar,\br,\gr $
are some constants and
\begin{eqnarray*}
b_0 & = & 2JH\w e_1 + JX_+\w e_- - JX_-\w e_+ \\
b_1 & = & X_+\w e_2 \\
b_2 & = & X_-\w e_+ .
\end{eqnarray*}
It is easy to see that $b_0= 2b_{e_2}$ (formula (\ref{bx})) and
$X_+e_2=0$, hence
$[b_0,b_0]=4\Om $  and
$$ [\ar b_0+\br b_1,\ar b_0+\br b_1] = \ar ^2 [b_0,b_0]=4\ar
^2\Om $$
(cf.~(\ref{bX})). Since
$$ [b_2,b_0]=2JX_-\w e_+\w e_1-2X_-\w e_+\w e_2,\qquad
[b_2,b_1]= -2H\w e_+\w e_2,\qquad [b_2,b_2]= 4X_-\w e_+\w e_1,$$
we have
$$[b,b]=\ar ^2 [b_0,b_0] + 2\gr\ar (2JX_-\w e_+\w e_1-2X_-\w e_+\w
e_2) -2\gr \br \cdot 2H\w e_+\w e_2 +4\gr ^2 X_-\w e_+\w e_1 .$$
The element
$$ 2[c,a] + ([b,b]-\ar ^2[b_0,b_0]) $$
is proportional to $\Om $ and has rank at most 4 (there are no
terms involving $JH$ and $ JX_+$), hence it is zero. In
particular (taking the term with $X_-\w e_+\w e_1$) we have $\gr
=0$. Finally we have
\be
b= \ar b_0 +\br b_1
\ee
and $[c,a]=0$. It follows that $Ha=0=X_+a$, hence $a=0$ (cf.~the
end of the previous section). This is the item 6 of the table.

\section{The proof for $c=0$}

We consider the case when $t=0$, hence equations
(\ref{cc})-(\ref{ab}) for $r=a+b$ reduce to
$$ [b,b]=0, \qquad [b,a]=0.$$
Since $b$ is a triangular $r$-matrix,
$$ b(\gotG^* ) = V_0 + \gotH _0,$$
where $\gotH _0:=b(V^*)\subset \gotH$, $V_0:=b(\gotH ^*)\subset
V$, is a Lie subalgebra of $\gotG = V\rtimes \gotH $. It follows
that $\gotH _0$  is a Lie subalgebra of $\gotH $ and $[\gotH _0
,V_0]\subset V_0$, therefore $b(\gotG ^*)=V_0\rtimes \gotH _0$.
Of course, $b$ is a triangular $r$-matrix on the smaller Lie algebra
$V_0\rtimes \gotH _0$. Let $b(\cdot )$ denote the linear bijection
from $V_0^*$ to $\gotH _0$ defined by $b$.
Equation $[b,b]=0$ is equivalent to
$$ [b(\xi ),b(\eta )] = b([\xi ,\eta ]_b), \qquad \xi,\eta\in V_0^*
$$
(cf.~(\ref{formula})). Applying the inverse map $f\colon \gotH
_0\to V_0^*$ of $b(\cdot )$ to the above equation changes it
from quadratic to a linear (!) one:
$$
f([X,Y]) = X f(Y) - Y f(X), \qquad X,Y\in \gotH _0,
$$
which says that $f$ is just a cocycle (on $\gotH _0$ with values
in $V_0^*$).

We consider four possible cases of $\dim V_0=\dim \gotH _0$
separately.

\subsection{$\dim V_0 =4$}

We shall show that there are no solutions of this type.  The
following lemma is not difficult.
\begin{Lem}
Any four-dimensional Lie subalgebra $\gotH _0$ of $\gotH
=sl(2,\bC )$ can be transformed by an internal automorphism to
\be\label{4dim}
 \left\{ \left( \ba{rr} z & w \\ 0 & -z \ea\right) : z,w\in
\bC \right\} = \left\langle H,JH,X_+,JX_+\right\rangle .
\ee
\end{Lem}
Assuming that $\gotH _0$ is given by (\ref{4dim}),
we are looking for cocycles $f\colon \gotH _0 \to V^*$. We
can replace $V^*$ by the isomorphic $\gotH $-module $V$. Set
$f(H)=:h$, $f(JH)=:k$, $f(X_+)=:x$ and $f(JX_+)=:y$.  The map
$f$ is a cocycle if and only if vectors $h,k,x,y$ satisfy
\begin{eqnarray}
Hk    & = & JH h   \nonumber \\
X_+ y & = & JX_+x  \nonumber \\
x     & = & Hx - X_+h \label{5.3} \\
y     & = & Hy - JX_+h \label{5.4} \\
y     & = & JHx - X_+k \nonumber \\
-x    & = & JHy - JX_+k \nonumber .
\end{eqnarray}
The first two equations are equivalent to $h=h^+e_+ + h^-e_-$,
$k=k^1e_1+k^2e_2$, $x= x^+e_+ + x^1e_1 + x^2 e_2$, $ y=y^+e_+
-x^2e_1 + y^2e_2$. Inserting this in (\ref{5.3}) gives $x^2=0$,
$x^1=-2h^-$. Inserting in (\ref{5.4}) gives $y^2=2h^-$. Then the
last two equations yield $y^+=-k^1$, $ x^+=k^2$. The general
solution is therefore as follows:
$$ h=h^+e_+ + h^-e_-,\;\;\;\; k=k^1e_1+k^2e_2,\;\;\;\;
x= k^2e_+ - 2h^-e_1,\;\;\;\; y=-k^1e_+ + 2h^-e_2 .
$$
These vectors are however linearly dependent:
$$ \det \left[\ba{cccc} h^-  &  0   &   0   &  0   \\
                       h^+  &  0   &  k^2 & - k^1  \\
		       0  &   k^1  & -2h^- & 0   \\
		       0  &  k^2   & 0   &   2 h^- \ea\right]=0,
$$
hence $f$ cannot be a bijection (this ends the proof).

\subsection{$\dim V_0 =3$}

There are three types of 3-dimensional subspaces $V_0$ of $V$:
\ben
\item {\em space-like}\, : \ $g|_{V_0}$ has signature $(0,3)$.
Then $\gotH _0\cong o(0,3)$.
\item {\em 3D-Minkowski}\, : \ $g|_{V_0}$ has signature $(1,2)$.
Then $\gotH _0\cong o(1,2)$.
\item {\em tangent to the light cone}\, : \ $g|_{V_0}$ has
signature $(0,2)$.
\een
In the first two cases $\gotH _0$ is simple and $f$ has to be
a coboundary:
$$f(X) = X\xi , \qquad X\in \gotH _0 \;\;\; (\mbox{for
some}\;\;\xi\in V_0^*).$$
Since each $\xi\in V^*$ has a nontrivial isotropy, $f$ cannot be
bijective.

In the third case we can assume the standard form $V_0=
\left\langle e_+, e_1,e_2\right\rangle $.
We have
$$ \gotH _0\subset \left\langle H,JH,X_+,JX_+\right\rangle ,
$$
because $\gotH _0$ is contained in the subalgebra stabilizing
$V_0$.
\begin{Lem} \ \
$\gotH _0 \supset \left\langle X_+,JX_+\right\rangle $.
\end{Lem}
\dowod
We set $\gotN := \left\langle X_+,JX_+\right\rangle $.
Since $\dim \gotH_0 =3$ and $\dim \gotN =2$, there exists $0\neq
Y\in \gotH _0\cap \gotN $. If $\gotH _0$ does not contain
$\gotN$, then $\gotH _0 +\gotN = \left\langle X_+,JX_+, \lr H +
\mu JH\right\rangle $, hence $JH\in \gotH_0 +\gotN $ and therefore
$$ JY = [JH,Y]\in \gotH _0$$
i.e. $\gotN = \left\langle Y,JY\right\rangle \subset \gotH _0$.

\qed

{}From the above lemma it follows that
\be\label{lam}
\gotH _0 = \left\langle X_+,JX_+, \lr H + \mu JH\right\rangle ,
\ee
where $\lr ^2 +\mu ^2 \neq 0$.
Let $(e^+,e^1,e^2)$ be the basis in $V_0^*$ dual to
$(e_+,e_1,e_2)$. The coordinates of an element $x\in V_0^*$ in
this basis are denoted by $x_+,x_1,x_2$.
We calculate also the action of $\gotH _0$ on $V_0^*$:
$$
\ba{rclrclrcl}
X_+e^+ & = & -e^1,\qquad & JX_+ e^+ & = & e^2,\qquad & (\lr H+\mu
JH)e^+ & = & -\lr e^+, \\
X_+e^1 & = & 0,\qquad & JX_+ e^1 & = & 0,\qquad & (\lr H+\mu
JH)e^1 & = & -\mu e^2 ,\\
X_+e^2 & = & 0,\qquad & JX_+ e^2 & = & 0,\qquad & (\lr H+\mu
JH)e^2 & = & \mu e^1.
\ea
$$
Let $f\colon \gotH _0\to V_0^*$ be a linear map and
$f(X_+)=:x$, $f(JX_+)=:y$, $f(\lr H +\mu JH)=: z$. It is a
cocycle if and only if
\begin{eqnarray}
X_+y & = & JX_+ x \nonumber \\
(\lr H +\mu JH)x - X_+z & = & \lr x + \mu y \label{pie}\\
(\lr H +\mu JH)y - JX_+z & = & -\mu x + \lr y .\label{dru}
\end{eqnarray}
The first equation is equivalent to $x_+=0=y_+$. Since $Hx=0=Hy$
and $X_+z=-z_+e^1$, $JX_+z=z_+e^2$, equations
(\ref{pie})--(\ref{dru}) are equivalent to
$$ \mu JH w +z_+(e^1-ie^2)=(\lr -i\mu )w,$$
where $w:= x+iy$ (just add (\ref{dru}) multiplied by $i$ to
(\ref{pie})), or to
\be\label{equa}
\mu (JH +i)w + z_+ (e^1-ie^2)=\lr w.
\ee
Since $JH (e^1-ie^2)= -i (e^1-ie^2)$, (\ref{equa}) is the
decomposition of $\lr w$ on components belonging to eigenspaces
of $JH$ (we know that $(JH)^2=-1$ on the subspace spanned by
$e^1,e^2$).
If $\lr =0$ then $z_+=0$ and $x,y,z$ are linearly dependent. In
order $f$ to be bijective we must have therefore $\lr \neq 0$.
In such a case we can assume in (\ref{lam}) and in the sequel
that $\lr =1$:
$$\mu (JH +i)w + z_+ (e^1-ie^2)= w.$$
Substituting here $w=w_{+i} +w_{-i}$, where $w_{+i}$ and
$w_{-i}$ are the eigenvectors of $JH$ corresponding to $+i$
and $-i$, respectively, we obtain $w_{+i}=0$. Therefore we have
$$ w = w_{-i} = z_+ (e^1-ie^2),$$
hence
$$
x=z_+e^1,\qquad y = -z_+e^2,\qquad z = z_+e^+ + z_1e^1+z_2e^2.
$$
Using the possibility of scaling $b$ (or $f$) by a non-zero
factor, we can assume that $z_+=1$:
\be
x=e^1,\qquad y = -e^2,\qquad z = e^+ + z_1e^1+z_2e^2.
\ee
Solving
$$ b_0(e^1)=X_+,\qquad b_0(-e^2)=JX_+,\qquad b_0 (e^+ +
z_1e^1+z_2e^2) = H +\mu JH,$$
we obtain
$$ b_0(e^1)=X_+,\qquad b_0(-e^2)=JX_+,\qquad b_0 (e^+) =
H+\mu JH -z_1X_+ + z_2JX_+,$$
hence finally
\be
b = e_1\wedge X_+ - e_2\wedge JX_+ + e_+ \wedge ( H+\mu JH -
z_1X_+ + z_2 JX_+).
\ee
Now note that
\be
b = b _{e_+} + e_+ \wedge ( \mu JH - z_1X_+ + z_2 JX_+).
\ee
Since $JH, X_+, JX_+$ belong to the isotropy subalgebra of
$e_+$, the above $b$ is of the form (\ref{bX}) and we can check
directly that $[b,b]=0$ (we know it already by the construction):
$$ [b,b] = [b_{e_+},b_{e_+}] = -g (e_+,e_+)=0.$$
We have two cases, depending on $\mu$:
\ben
\item $\mu\neq 0$. In this case one can get rid of $z_1,z_2$,
using the automorphisms generated by $X_+, JX_+$, since
$$X_+b= e_+\wedge \mu (-JX_+),\qquad JX_+ b = e_+\wedge\mu X_+
.$$
We have then $b=b_{e_+}+\mu e_+\wedge JH$. Since $JH$ has only
imaginary eigenvalues, by Prop.~\ref{noa}, adding $a$ does not
lead to new solutions, hence we get item 7 of the table.
\item $\mu = 0$. The one-parameter group of automorphisms generated
by $JH$ acts on $b$ according to the linear system of
differential equations $\dot{z_1} = z_2$, $\dot{z_2} =-z_1$.
Therefore we can assume that $z_2=0$: $b= b_{e_+} + ze_+\wedge
X_+$. Again, there is no need to consider nontrivial $a$, since
$X_+$ is nilpotent. We get then item 8 of the table.
\een

\subsection{$\dim V_0 =2$}

There are three normal forms of a 2-dimensional subspace $V_0$
of $V$:
\ben
\item $V_0=\left\langle e_1,e_2 \right\rangle $ \ ({\em
space-like}\, : \ $g|_{V_0}$ has signature $(0,2)$).
Then $\gotH _0 =\left\langle H,JH \right\rangle $.
\item $V_0=\left\langle e_+,e_- \right\rangle $ \ ({\em
2D-Minkowski}\, : \ $g|_{V_0}$ has signature $(1,1)$).
Then $\gotH _0=\left\langle H,JH \right\rangle $.
\item $V_0=\left\langle e_1,e_+ \right\rangle $ \ ({\em
tangent to the light cone}\, : \ $g|_{V_0}$ has signature $(0,1)$).
Then $\gotH _0\subset \left\langle H,X_+,JX_+ \right\rangle $.
\een
(The simplest way to prove it is to note that 2-dimensional
subspaces of $V$ correspond to simple bivectors, i.e. some
elements of $\gotH$; the classification of the latter is easy.)

In the first case, $b=x\wedge H+y\wedge JH$, where $x,y\in
\left\langle e_1,e_2 \right\rangle $. We have
$$ \frac12 [b,b] = y\wedge JH y\wedge JH + y\wedge JH x\wedge H,$$
hence $[b,b]=0$ implies the linear dependence of $y,JHy$, i.e.
$y=0$. This is in contradiction with $\dim V_0=2$.

In the second case, $b=x\wedge H+y\wedge JH$, where $x,y\in
\left\langle e_+,e_- \right\rangle $. We have
$$ \frac12 [b,b] = x\wedge H x\wedge H + x\wedge H y\wedge JH,$$
hence $[b,b]=0$ implies $x\wedge Hx=0=x\wedge Hy$. Since we
consider only nonzero $x, y$, this means that there exist
$\lr ,\mu$ such that $x=\lr Hx$ and $x=\mu Hy$. We have
therefore $x=\lr \mu H^2 y =\lr\mu y$. This is in contradiction
with $\dim V_0=2$.

In the third case, $b$ is of the following form:
$$ b= x\wedge X_+ + y \wedge JX_+ + z \wedge H ,$$
where $x=x^+ e_+ + x^1e_1$, etc. A simple calculation shows that
$[b,b]=0$ if and only if
\begin{eqnarray*}
x^1 (x^1- z^+) + 2x^+z^1 & = & 0 \\
y^1 (x^1- z^+) + 2y^+z^1 & = & 0 \\
z^1 (x^1- z^+) + 2z^+z^1 & = & 0 .
\end{eqnarray*}
Note that if
$$\left[\ba{c} x^1- z^+ \\ 2z^1 \ea\right]$$
is a non-zero vector, then $x,y,z$ are in the same
one-dimensional subspace. This would mean that $\dim V_0\leq 1$.
We conclude that $x^1- z^+ =0 =z^1$ and
\be
b= e_1\wedge (x^1X_++y^1JX_+) + e_+ \wedge ( x^+X_+ + y^+JX_+ +
x^1H) .
\ee
Now we shall reduce the number of parameters, acting by suitable
automorphisms. We consider separately two cases.

\vspace{1mm}

\noindent
{\bf Case 1.}\ \ $x^1\neq 0$. \

Since  $JX_+\, b = - x^1 e_+\wedge JX_+$ (which means
$\dot{y}^+=-x^1=const$) and $x^1\neq 0$, we can pass to the
situation when $y^+=0$. Using another group of automorphisms,
the one generated by $H$, we get the change of parameters as
follows
$$ \dot{x}^1 = x^1,\qquad \dot{y}^1=y^1,\qquad \dot{x}^+=2x^+.$$
Using this and the possibility of multiplying $b$ by a nonzero
number, we get
\be
b= e_1\wedge (X_++y^1JX_+) + e_+ \wedge (H + x^+X_+),
\ee
where $x^+=0,\pm 1$. For $v\in V$ we have
$$ (-v)b = (v^1-y^1v^2-2x^+v^-)e_1\wedge e_+ -2y^1v^-e_1\wedge
e_2 + v^-e_-\wedge e_+ ,$$
hence we can assume that $a$ is of the form
$$ a = \ar e_+\wedge e_2 + e_-\wedge (\gr _1 e_1+\gr _2 e_2) +
\mu e_1\wedge e_2 $$
(no component with $e_-\wedge e_+$, $e_1\wedge e_+$). A simple
calculation yields
$$ [a,b] = (2\gr _1 - y^1\gr _2)e_1\wedge e_-\wedge e_+ +
(\mu -2x^+\gr _2)e_1\wedge e_+\wedge e_2 -\gr _2 e_+\wedge
e_-\wedge e_2 .$$
It follows that $[a,b]=0$ if and only if $a= \ar e_+\wedge e_2$,
which is item 9 of the table.

\vspace{1mm}

\noindent
{\bf Case 2.}\ \ $x^1= 0$. \

In this case we have
\be
b= y^1 e_1\wedge JX_+ + e_+ \wedge ( x^+X_+ + y^+JX_+ )
\ee
$y^1\neq 0\neq x^+$ (because $\dim V_0 =2$).
Since  $X_+\, b = - y^1 e_+\wedge JX_+$ (which means
$\dot{y}^+=-y^1=const$) and $y^1\neq 0$, we can pass to the
situation when $y^+=0$. Using another group of automorphisms,
the one generated by $H$, we get the change of parameters as
follows
$$ \dot{y}^1 = y^1,\qquad \dot{x}^+=2x^+.$$
Using this and the possibility of multiplying $b$ by a nonzero
number, we get
\be\label{pmb}
b=  \pm e_1\wedge JX_+ + e_+ \wedge X_+ .
\ee
Now, observe that the reflection $e_2\mapsto -e_2$ (other
elements of the basis unchanged) yields an automorphism of
$\gotG$ which on $\gotH$ coincides with the `complex
conjugation' (if the real part is spanned by $H,X_+,X_-$), in
particular $JX_+\mapsto -JX_+$. It means that we can choose plus
sign in (\ref{pmb}):
\be
b=  e_1\wedge JX_+ + e_+ \wedge X_+.
\ee
 For $v\in V$ we have
$$ (-v)b = (v^2-2v^-)e_1\wedge e_+ -2v^-e_1\wedge e_2 ,$$
hence we can assume that $a$ is of the form
$$ a = e_-\wedge (\ar _1 e_1+\ar _2 e_2) +\gr e_-\wedge e_+ +
\ar e_+\wedge e_2 $$
(no component with $e_1\wedge e_+$, $e_1\wedge e_2$). A simple
calculation yields
$$ [a,b] = - \ar _2 e_1\wedge e_-\wedge e_+ +
(2\gr - 2\ar _2)e_1\wedge e_+\wedge e_2 .$$
It follows that $[a,b]=0$ if and only if $a= \ar _1e_-\wedge e_1
+ \ar  e_+\wedge e_2$,
which is item 10 of the table.

\subsection{$\dim V_0 =1$}

In this case $b=v\wedge X$ for some nonzero $v\in V$, $X\in
\gotH$. Since $X$ has to preserve $V_0:=\left\langle
v\right\rangle$, $v$ is an eigenvector of $X$ and $[b,b]=0$
automatically in this case. We can always rescale $X$ in such a
way that $Xv=0$ or $Xv=v$.

The classification procedure is simple. Any nilpotent $X$ is
equivalent to $X_+$ and any semisimple $X$ is equivalent to $\lr
H+\mu JH$. We have then the following possibilities:
$$
\ba{c|c}
X & v \\
\hline
X_+ & v\in \left\langle e_+,e_2\right\rangle \\
JH &  v\in \left\langle e_+,e_-\right\rangle \\
H  & v\in \left\langle e_1,e_2\right\rangle,\; v=e _{\pm } \\
H+\br JH \;\; \br\neq 0 & v=e _{\pm }.
\ea
$$
 Note that we can still restrict the
possibilities. Namely, we use the automorphisms generated by
$JX_+$, $H$, $JH$ (and scaling) in cases when $X=X_+$, $X=JH$,
$X=H$, respectively, to pass from two-dimensional
eigenspaces of $X$ to specific vectors: $e_+,e_2$ in the first
case, $e_{\pm},e_0,e_3$ in the second case and $e_{\pm},e_1$ in
the third. We also use the reflection $e_3\mapsto -e_3$ in order
to replace $e_-\w JH$, $e_-\w H$, $e_-\w (H+\br JH)$ by
$e_+\w JH$, $e_+\w H$, $e_+\w (H-\br JH)$, respectively.

The results are presented in the following table, where we have
also shown which $a$ satisfy
$$ [a,b]= v\wedge Xa =0,$$
how they can be simplified using $(-v)b$ and which still can be
simplified using $H$ (in one case also $JH$) to get the final
number of parameters $\# $. This covers items 11--18 in
Table~\ref{tabela}.
\pagebreak
\begin{table}[t]
\bt{|c|c|c|c|c|}
\hline
$b$  &  $a$ belongs to  &  $a+(-v)b $ & still use & \#  \\
\hline
$e_2\w X_+$ & $\lel e_+\w e_1, e_+\w e_2, e_-\w e_2, e_1\w e_2\rr
$ & $\lel e_+\w e_1, e_-\w e_2\rr $ & $H$ & 1 \\
$e_+\w X_+$ & $\lel e_-\w e_+, e_{\pm}\w e_1, e_{\pm}\w e_1\rr $
&  $\lel e_-\w e_+, e_-\w e_1, e_{\pm}\w e_2\rr $ & $H$ & 3 \\
$e_0\w JH$ & $\lel e_0\w e_1,e_0\w e_2,e_0\w e_3, e_1\w e_2\rr $
& $\lel e_0\w e_3, e_1\w e_2\rr $ &  & 2 \\
$e_3\w JH$ & $\lel e_0\w e_3,e_1\w e_3,e_2\w e_3, e_1\w e_2\rr $
& $\lel e_0\w e_3, e_1\w e_2\rr $ &  & 2 \\
$e_+\w JH$ & $\lel e_+\w e_1,e_+\w e_2,e_+\w e_3, e_1\w e_2\rr $
& $\lel e_0\w e_3, e_1\w e_2\rr $ & $H$ & 1 \\
$e_1\w H$ & $\lel e_-\w e_+,e_{\pm}\w e_1, e_1\w e_2\rr $
& $\lel e_0\w e_3, e_1\w e_2\rr $ &  & 2 \\
$e_+\w H $ & $\lel e_+\w e_1,e_+\w e_2,e_+\w e_3, e_1\w e_2\rr $
& $\lel  e_+\w e_1,e_+\w e_2,e_1\w e_2\rr $ & $JH$, $H$ & 1 \\
$e_+\w (H+\br JH) $
& $\lel e_+\w e_1,e_+\w e_2,e_+\w e_3, e_1\w e_2\rr $
& $\lel  e_1\w e_2\rr $ & $H$ & 1 \\
\hline
\et
\caption{The lowest non-zero rank of $b$}
\end{table}

\subsection{$b_0 =0$}

The classification of $a\in \Ldwa V$ is the same as the
classification of elements of $sl (2,\bC )$. Additionally, we
identify proportional elements.
 The normal forms are $X_+\sim e_1\w e_+$, $JH\sim e_1\w e_2$ and
 $H+\ar JH\sim e_0\w e_3 +\ar e_1\w e_2$ (items 19--21).

\section{Final remarks}

\ben
\item Unfortunately, we were not able to solve generally the
`classical modified Yang-Baxter equation' $[b,b]=t\Om ,\; t\neq
0$, in spite of the existence of general solution in the case of
simple Lie algebras given by Belavin and Drinfeld \cite{BD}.
\item According to Remark~1.8~of \cite{qpoi}, any solution of
(\ref{bb}), (\ref{ab}) with $c=0$ (non-deformed classical
Lorentz subgroup) defines directly a quantum
Poincar\'{e} group. All non-zero solutions with $c=0$ and $t=0$
are given as items 7--21 of Table~\ref{tabela}. Some solutions
with $c=0$ and $t\neq 0$ are given in (\ref{be0+}), (\ref{be1+}),
(\ref{be1+a}).
\item Poisson structures on the Poincar\'{e} group acting in
2-dimensional space-time have been classified in \cite{poi}.
\item The 3-dimensional case is investigated in \cite{stach}.
\item For each Poisson Poincar\'{e} group $G$ there is exactly one
Poisson Minkowski space $M$ (with a Poisson action of $G$ on
$M$), cf.~\cite{poihom,standr}.
\item Some classical-mechanical models of particles based on
Poisson Poincar\'{e} symmetry were discussed in
\cite{poi,poican,k-part}. For a short review see \cite{zakop}.
\een

\section*{Acknowledgments}

The author would like to thank to Dr P.~Podle\'{s} and
Dr~F.~Burstall for valuable discussions.

This research was supported by Polish KBN grant No. 2 P301 020 07.


\begin{thebibliography}{99}

\bibitem{PPgr} S.~Zakrzewski, {\em Poisson Poincar\'{e} groups},
in: ``Quantum Groups, Formalism and Applications'',
Proceedings of the XXX Winter School on Theoretical Physics
14--26 February 1994, Karpacz, J.~Lukierski, Z.~Popowicz,
J.~Sobczyk (eds.), Polish Scientific Publishers PWN, Warsaw
1995, pp.~433--439 (also: hep-th{/}9412099).

\bibitem{qpoi} P. Podle\'{s} and S.L. Woronowicz, {\em On the
classification of quantum Poincar\'{e} groups},
hep-th{/}9412059, UC Berkeley preprint, PAM-632, to appear in
Comm.~Math.~Phys.

\bibitem{ccc} W. Greub, S. Halperin and R. Van Stone, {\em
Curvature, Connection and Cohomology}, Pure and Applied
Mathematics, vol. 47, III, Academic Press, New York, 1976.

\bibitem{D:ham} V. G. Drinfeld, {\em Hamiltonian structures on
Lie groups, Lie bialgebras and the meaning of the classical
Yang-Baxter equations},  Soviet Math. Dokl. {\bf 27} (1983),
68--71.

\bibitem{D} V. G. Drinfeld, {\em Quantum groups}, Proc. ICM,
Berkeley, 1986, vol.1, 789--820.

\bibitem{S-T-S} M. A. Semenov-Tian-Shansky, {\em Dressing
transformations and Poisson Lie group actions},
Publ. Res. Inst. Math. Sci.,
Kyoto University {\bf 21} (1985), 1237--1260.

\bibitem{Lu-We} J.-H. Lu and A. Weinstein, {\em Poisson Lie Groups,
Dressing Transformations and Bruhat Decompositions},
J. Diff. Geom. {\bf 31} (1990), 501--526.

\bibitem{Bourb} N. Bourbaki, {\em Groupes et alg\`{e}bres de
Lie}, Chapitre VIII, Hermann 1975.

\bibitem{inho}  P. Podle\'{s} and S.L. Woronowicz, {\em On the
structure of inhomogeneous quantum groups},
hep-th{/}9412058, UC Berkeley preprint, PAM-631.

\bibitem{abel} S. Zakrzewski, {\em Geometric quantization of
Poisson groups --- diagonal and soft deformations},
Proceedings of the Taniguchi Symposium {\sl
Symplectic geometry and quantization problems}, Sanda (1993),
Y.~Maeda, H.~Omori and A.~Weinstein (Eds.), Contemporary
Mathematics {\bf 179}, 1994, 271--285.

\bibitem{repoi} S. Zakrzewski,
{\em Poisson structures on the Lorentz group},
 Lett. Math. Phys. {\bf 32} (1994), 11--23.

\bibitem{luki} S. Zakrzewski,
{\em Quantum Poincar\'{e} group related to $\kappa$-Poincar\'{e}
algebra},  J. Phys. A: Math. Gen. {\bf 27} (1994), 2075--2082.

\bibitem{BD} A. Belavin and V.G. Drinfeld, {\em Triangle
equations and simple Lie algebras}, Sov.~Sci.~Rev.~Math.~{\bf 4}
 (1984), 93--165.

\bibitem{poi} S. Zakrzewski,
{\em Poisson space-time symmetry and corresponding
elementary systems},
in: ``Quantum Symmetries'', Proceedings of the II International
Wigner Symposium, Goslar 1991, H.D.~Doebner and V.K.~Dobrev (Eds.),
pp. 111--123.

\bibitem{stach} P.~Stachura, in preparation.

\bibitem{poihom} S.~Zakrzewski, {\em Poisson homogeneous spaces},
in: ``Quantum Groups, Formalism and Applications'',
Proceedings of the XXX Winter School on Theoretical Physics
14--26 February 1994, Karpacz, J.~Lukierski, Z.~Popowicz,
J.~Sobczyk (eds.), Polish Scientific Publishers PWN, Warsaw
1995, pp.~629--639 (also hep-th{/}9412101).

\bibitem{standr} S.~Zakrzewski, {\em Phase spaces associated
with standard $r$-matrices}, q-alg{/}9511002, Warsaw preprint
1995.

\bibitem{poican} S. Zakrzewski, {\em Poisson Poincar\'{e}
particle and canonical variables}, in: ``Generalized
Symmetries'',
Proceedings of the International Symposium on Mathematical
Physics, Clausthal, July 27--29, 1993, H.-D. Doebner, V.K.
Dobrev and A.G. Ushveridze (Eds.), 1994, pp. 165--171.

\bibitem{k-part} S.~Zakrzewski, {\em On the classical $\kappa
$-particle},
in: ``Quantum Groups, Formalism and Applications'',
Proceedings of the XXX Winter School on Theoretical Physics
14--26 February 1994, Karpacz, J.~Lukierski, Z.~Popowicz,
J.~Sobczyk (eds.), Polish Scientific Publishers PWN, Warsaw
1995, pp.~573--577 (also: hep-th{/}9412098).

\bibitem{zakop} S.~Zakrzewski, {\em Classical mechanical systems
based on Poisson symmetry}, submitted for the Proceedings of
The Second German--Polish Symposium
``New Ideas in the Theory of Fundamental Interactions'',
September 1995, Zakopane, Poland.

\end{thebibliography}
\end{document}